\documentclass{article}

\usepackage[latin1]{inputenc}
\usepackage[T1]{fontenc}
\usepackage{textcomp}
\usepackage{graphicx}
\usepackage{epsfig}
\usepackage{amsmath}
\usepackage{amssymb}
\usepackage{amsthm}
\usepackage{booktabs}
\usepackage{stmaryrd}
\usepackage{url}
\usepackage{longtable}
\usepackage[figuresright]{rotating}
\usepackage{bbm}
\usepackage[all]{xy}
\usepackage{pdfsync}
\usepackage[toc,page]{appendix}

\usepackage{hyperref}

\usepackage{makeidx}
\makeindex

\newcommand{\komplex}{\mathbb{C}}
\newcommand{\complex}{\mathbb{C}}
\newcommand{\ganz}{\mathbbm{Z}}
\newcommand{\reell}{\mathbbm{R}}

\newcommand{\Q}{\mathbbm Q}
\newcommand{\q}{\mathbbm Q}
\newcommand{\qbar}{{\overline{\mathbb Q}}}
\newcommand{\diagram}[1]{\begin{equation}
\xymatrix{#1}
\end{equation}}

\newcommand{\diagramlab}[2]{\begin{equation}
\label{#1}
\xymatrix{#2}
\end{equation}}

\newcommand{\diagramst}[1]{\begin{equation*}
\xymatrix{#1}
\end{equation*}}

\newcommand{\adfin}{\mathbb {A}_{f} }
\newcommand{\real}{\mathbbm R}

\newcommand{\matrt}[1]{  \left( \begin{array}{cc} #1 \end{array} \right)}
\newcommand{\sieg}{\mathbb H _g}
\newcommand{\eq}[1]{\begin{align}#1\end{align}}
\newcommand{\eqst}[1]{\begin{align*}#1\end{align*}}
\newcommand{\siegel}{{Sh}}
\newcommand{\al}[1]{\begin{align}#1\end{align}}
\newcommand{\alst}[1]{\begin{align*}#1\end{align*}}

\theoremstyle{plain}
\newtheorem{theorem}{Theorem}[section]

\newtheorem{criterion}[theorem]{Criterion}
\newtheorem{proposition}[theorem]{Proposition}

\newtheorem{lemma}[theorem]{Lemma}

\theoremstyle{remark}
\newtheorem{remark}{Remark}

\theoremstyle{definition}
\newtheorem{definition}{Definition}[section]

\newenvironment{keywords}
   {\begin{trivlist}\item[]{Key words:}\ }
   {\end{trivlist}}

\newenvironment{msc2010}
   {\begin{trivlist}\item[]{2010 Mathematics Subject Classification:}\ }
   {\end{trivlist}}

\title{On Bost-Connes type systems and Complex Multiplication}
\author{Bora Yalkinoglu
}
\date{}

\begin{document}
\maketitle

\begin{abstract}
By using the theory of Complex Multiplication for general Siegel modular varieties we construct arithmetic subalgebras for BC-type systems attached to number fields containing a CM field. The abelian extensions obtained in this way are characterized by results of \cite{wei}. Our approach is based on a general construction of BC-type systems of Ha and Paugam \cite{HaPa} and extends the construction of the arithmetic subalgebra of Connes, Marcolli and Ramachandran \cite{cmr1} for imaginary quadratic fields. 
\end{abstract}

\setcounter{tocdepth}{1}
\tableofcontents

\begin{msc2010}
58B34, 11R37, 14G35
\end{msc2010}


%%%%%%%%%%%%%%%%%%%%%%%%%%%%%%%%%%%%%%%%%%%%%
%%%%%%%%        INTRODUCTION
%%%%%%%%%%%%%%%%%%%%%%%%%%%%%%%%%%%%%%%%%%%%%
\section{Introduction}

In their fundamental article \cite{bc} Bost and Connes constructed a quantum statistical mechanical system, the so called Bost-Connes or BC system, that recovers (among other arithmetic properties of $\q$) the explicit class field theory of $\q$. It is a natural question to ask for other quantum statistical mechanical systems that recover - at least partially - the (explicit) class field theory of number fields different from $\q$. \\
This paper is a contribution to this question in the case of CM fields. In nature a CM field arises as (rational) endomorphism ring of an Abelian variety with Complex Multiplication. 
\subsection*{Background}
\label{background}
A quantum statistical mechanical system $\mathcal A = (A,(\sigma_t)_{t\in \real})$ is a $C^*$-algebra $A$ together with a one-parameter group of automorphisms $(\sigma_t)_{t\in \real}$. Sometimes we call $\mathcal A$ simply a $C^*$-dynamical system. One should think of $A$ as the algebra of observables of a quantum physical system with the one-parameter group $(\sigma_t)_{t\in \real}$ implementing the time evolution of the observables. A symmetry of $\mathcal A$ is a $C^*$-automorphism of $A$ which commutes with the time evolution.
Interesting properties of $\mathcal A$ can be read off its equilibrium states at a temperature $T \in [0,\infty)$. An equilibrium state at inverse temperature $\beta = \frac 1 T$ is given by a so called $KMS_\beta$-state. See \cite{brat1} and \cite{brat2} for more information. \\  \\
Now a first step towards generalizing the Bost-Connes system to other number fields is given by a system of the following kind (see \cite{cmr2} and \cite{Sergey}) \\ \\
A \textbf{BC-type system} for a number field $K$ is a quantum statistical mechanical system $\mathcal A = (A,(\sigma_t)_{t\in \real})$ with the following properties 
\begin{itemize}
\item[(i)] The partition function of $\mathcal A$ is given by the Dedekind zeta function of $K$. 
\item[(ii)] The quotient of the idele class group $C_K$ by the connected component $D_K$ of the identity of $C_K$ acts as symmetries on $\mathcal A$. 
\item[(iii)] For each inverse temperature $0 < \beta \leq 1$ there is a unique $KMS_\beta$-state.
\item[(iv)] For each $\beta > 1$ the action of the symmetry group $C_K / D_K$ on the set of extremal $KMS_\beta$-states is free and transitive.
\end{itemize}
\begin{remark}
i) Using class field theory we see that property $(2)$ simply states that the Galois group $Gal(K^{ab} / K) \cong C_K / D_K$ of the maximal abelian extension $K^{ab}$ of $K$ is acting by symmetries on $\mathcal A$. \\ 
ii) Properties $(3)$ and $(4)$ say that there is a spontaneous symmetry breaking phenomenon at $\beta = 1$ (see \cite{bc} or pp. 400 in the book \cite{conmar} of Connes and Marcolli).
\end{remark}
Thanks to the work of Ha and Paugam \cite{HaPa} BC-type systems are known to exist for an arbitrary number field $K$, we denote their solution (see the next section for the definition) by \eq{\label{bchapa1} \mathcal A _K = (A_K,(\sigma_t)_{t \in \real}).} 
The paper of Laca, Larsen and Neshveyev \cite{Sergey} gives a different description of $\mathcal A _K$. In fact \cite{HaPa} yields that $\mathcal A _K$ fulfills the first two properties and \cite{Sergey} verifies the last two properties. \\
We should mention that the approach of \cite{HaPa} works much more general, it allows to attach a quantum statistical mechanical system to an arbitrary Shimura variety; the above construction being a special case.
The Shimura-theoretic approach of \cite{HaPa} is based on the foundational work of Connes, Marcolli and Ramachandran (see \cite{cmr1} and \cite{cmr2}). \\ \\
The notion of an \textbf{arithmetic subalgebra} for a BC-type system, first encountered in \cite{bc}, hints to a possible relation between seemingly unrelated areas, namely it asks for a connection between quantum statistical mechanical systems (physics/operator algebras)  and class field theory (number theory). More precisely (see \cite{cmr2}) \\ \\
A \textbf{BC-system} for a number field $K$ is a BC-type system $\mathcal A = (A,(\sigma_t)_{t\in \real})$ such that
\begin{itemize}
\item[(v)] There is a $K$-rational subalgebra $A^{arith}$ of $A$, called arithmetic subalgebra of $\mathcal A$, such that for every extremal $KMS_\infty$-state $\varrho$ and every $f \in A^{arith}$ we have \eq{\label{arithalgdef}\varrho(f) \in K^{ab}} and further $K^{ab}$ is generated over $K$ in this manner, i.e. \eq{K^{ab}=K(\varrho(f) \ | \ \text{$\varrho$ extremal $KMS_\infty$-state, }f \in A^{arith} ).}
\item[(vi)] Let $\nu \in C_K$ be a symmetry of $\mathcal A$, $\varrho$ an extremal $KMS_\infty$-state of $\mathcal A$ and $f \in A^{arith}$ and denote further by $[\nu] \in Gal(K^{ab}/K)$ the image of $\nu$ under Artin's reciprocity morphism $C_K \to Gal(K^{ab}/K)$. If we denote the action (given by pull-back) of $\nu$ on $\varrho$ by ${}^{\nu}\hspace{-0.6mm}\varrho$, we have the following compatibility relation \eq{\label{natactidel}{}^{\nu}\hspace{-0.6mm}\varrho(f)  = [\nu]^{-1}(\varrho(f)) \in K^{ab}.}
\end{itemize}

Now the difficulty of constructing BC-systems comes from its relation with Hilbert's 12th problem,
which asks for an explicit class field theory of a number field $K$. This problem is completely solved only in two cases, namely the case of $K = \q$ and the case of $K$ equal to an imaginary quadratic field, like $K = \q(i)$.
It does not come as a surprise that BC-systems are so far known to exist only in these two cases. For $K=\q$ see \cite{bc} and for the imaginary quadratic case we refer to \cite{cmr1} and the very detailed exposition given in \cite{conmar} pp. 551.  \\ \\
The construction in  \cite{cmr1} is based on the theory of Complex Multiplication (see section \ref{cmtheory}) which is a part of the arithmetic theory of Shimura Varieties. This theory allows to construct explicitly abelian extensions of higher dimensional generalizations of imaginary quadratic number fields, so called CM fields. A CM field $E$ is a totally imaginary quadratic extension of a totally real number field. So for example cyclotomic number fields $\q(\zeta_n)$, $\zeta_n$ being a primitive n-th root of unity, are seen to be CM-fields. \\
Abelian extensions of a CM field $E$ are obtained by evaluating arithmetic Modular functions on so called CM-points on a Siegel upper half plane (see section \ref{intro2} for more information). \\
Except for the case of an imaginary quadratic field, it is unfortunately not possible to generate the maximal abelian extension $E^{ab}$ of a CM field $E$ in this way, but still
a non trivial abelian extension of infinite degree over $E$. We denote the latter abelian extension of $E$ by \eq{\label{extcm}E^c \subset E^{ab}.} See Thm. \ref{bigtheorem} for the characterization of $E^c$ given in \cite{wei}. As we remarked above for $E$ equal to an imaginary quadratic field there is an equality $E^c=E^{ab}$. \\ \\
Due to the lack of a general knowledge of the explicit class field theory of an arbitrary number field $K$ we content ourself with the following weaking of a BC-system: \\ \\
Let $F \subset K^{ab}$ be an arbitrary abelian extension of $K$. A \textbf{partial BC-system} for the extension $F / K$ is defined like a BC system $\mathcal A = (A,(\sigma_t)_{t\in \real})$ for $K$ except that we do not demand the arithmetic subalgebra $A^{arith}$ of $A$ to generate the maximal abelian extension $K^{ab}$ of $K$ but instead the abelian extension of $F$, i.e. in property $(v)$ we replace $K^{ab}$ by $F$ and call the corresponding $K$-rational subalgebra $A^{arith}$ a (partial) arithmetic subalgebra of $\mathcal A$.

\subsection*{Statement of our result}
Now the aim of our paper is to prove the following 
\begin{theorem}
\label{theorem}
Let $K$ be a number field containing a CM field. Denote by $E$ the maximal CM field contained in $K$ and define the abelian extension $K^c$ of $K$ to be the compositum  \eq{\label{kc}K^c = K \cdot E^c} (with the notation from (\ref{extcm})). \\
Then the BC-type system $\mathcal A_K$ of \cite{HaPa} (see (\ref{mathcalAK})) is a partial BC-system for the extension $K^c / K$, i.e. there exists a (partial) arithmetic subalgebra $A_K^{arith}$ for $\mathcal A _K$.
\end{theorem}
We want to explain our construction of $A^{arith}_K$ (cf. section \ref{constrarith}). It is inspired by the construction given in \cite{cmr1} and \cite{cmr2} (see also pp. 551 \cite{conmar}). 
\subsection*{Idea of our construction}
Let $(G,X,h)$ be a Shimura datum and denote by $Sh(G,X,h)$ the associated Shimura variety.
In \cite{HaPa} the authors associate to the datum $(G,X,h)$, the variety $Sh(G,X,h)$ and some additional data a quotient map \eq{U \longrightarrow Z} between a topological groupoid $U$ and a quotient $Z = \Gamma \backslash U$ for a group $\Gamma$. Out of this data they construct a $C^*$-dynamical system $\mathcal A = (A,(\sigma_t)_{t\in \real})$. A dense $*$-subalgebra $H$ of $A$ is thereby given by the compactly supported, continuous functions \eq{H = C_c(Z)} on $Z$, where the groupoid structure of $U$ induces the $*$-algebra structure on $H$. \\ \\
Moreover there are two variations of the above quotient map which can be put into the following (commutative) diagram
\diagram{U \ar[rr] \ar @{} [drr] | { } & & Z \\ U^+ \ar[u] \ar[rr] \ar[d] && Z^+ \ar[u]_{} \ar[d] \\ U^{ad} \ar[rr] && Z^{ad}.}
For more information see section \ref{hapabcm}. We will apply this general procedure in two special cases. For this we fix a number field $K$ together with a maximal CM subfield $E$. 
\subsubsection*{I) The BC-type system $\mathcal A_K$} 
The (0-dimensional) Shimura datum $\mathcal S _K = (T^K,X_K,h_K)$ gives rise to a quotient map denoted  by \eqst{U_K \longrightarrow Z_K.} In this case the groupoid $U_K$ is of the form (cf. (\ref{notlln})) \eqst{U_K = T^K(\adfin) \boxplus (\widehat {\mathcal O}_K \times Sh(T^K,X_K,h_K)).}
The associated $C^*$-dynamical system is denoted by \eq{\mathcal A _K =(A_K,(\sigma_t)_{r\in \real}).} It gives rise to a BC-type system for $K$. For the precise definition of $\mathcal A _K$ and its properties (e.g. symmetries, extremal $KMS_\infty$-states) we refer the reader to section \ref{shdatum1} and \ref{mathcalAK}. Moreover we denote by $H_K$ the dense subalgebra of $A_K$ given by \eq{H_K = C_c(Z_K).}
\subsubsection*{II) The Shimura system $\mathcal A_\siegel$}
To the CM field $E$ we associate the Shimura datum \eqst{\mathcal S_\siegel = (GSp(V_E,\psi_E),\mathbb H^±_g,h_{cm})} where in fact the construction of the morphism $h_{cm}$ takes some time (see \ref{construction}). This is due to two difficulties that arise in the case of a general CM field $E$ which are not visible in the case of imaginary quadratic fields. On the one hand one has to use the Serre group $S^E$ and on the other hand in general the reflex field $E^*$ of a CM field $E$ is not anymore equal to $E$ (see \ref{reflexB}). We denote the associated quotient map by \eqst{U_\siegel \longrightarrow Z_\siegel} and analogously its variations (see \ref{somequotientmaps}). Here the relevant groupoids in hand are of the following form \eqst{U_\siegel = GSp(\adfin) \boxplus (\Gamma_{\siegel,M} \times Sh(GSp(V_E,\psi_E),\mathbb H^±_g,h_{cm}) ) } and \eqst{U_\siegel^{ad}=GSp^{ad}(\q)^+ \boxplus (\Gamma^{ad}_{\siegel,M} \times \mathbb H_g). } 
We denote the resulting $C^*$-dynamical system by $\mathcal A _\siegel$ and call it Shimura system (cf. section \ref{bcmalgebras}).
\begin{remark}
In the case of an imaginary quadratic field $K$ the Shimura system $\mathcal A_\siegel$  gives rise to the $GL_2$-system of Connes and Marcolli \cite{conmar}.
\end{remark}
The second system is of great importance for us because of the following: Denote by $x_{cm} \in \mathbb H_g$ the CM-point associated with $h_{cm}$ and denote by $\mathcal M^{cm}$ the ring of arithmetic Modular functions on $\mathbb H_g$ defined at $x_{cm}$. By the theory of Complex Multiplication we know that for every $f \in \mathcal M^{cm}$ we have (cf. (\ref{kc}) and \ref{defmathcalM}) \eqst{f(x_{cm})\in K^{c}\subset K^{ab}} and moreover $K^c$ is generated in this way. Our idea is now that $\mathcal M ^{cm}$ gives rise to the arithmetic subalgebra $A^{arith}_K$. More precisely we will construct a (commutative) diagram (see \ref{thebcmmap}) \diagram{U_K \ar[rr] \ar[d] && Z_K \ar[d] \\ U_\siegel \ar[rr] && Z_\siegel} which is induced by a morphism of Shimura data $\mathcal S _K \to \mathcal S _\siegel$ constructed in section \ref{themapphi}. \\
Then, using criterion \ref{crit}, we see that the morphism $Z^+_\siegel \to Z_\siegel$ (from the above diagram) is invertible and obtain in this way a continuous map \eqst{\Theta : Z_K \longrightarrow Z_\siegel \longrightarrow Z^+_\siegel \longrightarrow Z_\siegel^{ad}.}
Using easy properties of the automorphism group of $\mathcal M^{cm}$, we can model each $f$ in $\mathcal M^{cm}$ as a function $\widetilde f$ on the space $Z^{ad}_\siegel$ (which might have singularities). Nevertheless in Prop. \ref{pullbackprop} we see that for every $f \in \mathcal M^{cm}$ the pull back $\widetilde f \circ \Theta$ lies in $H_K = C_c(Z_K)$ and we can define $A^{arith}_K$ as the $K$-algebra generated by these elements, i.e. \eqst{A^{arith}_K = \ <\widetilde f \circ \Theta \ | \ f \in \mathcal M^{cm} >_K.}
Now, using the classification of extremal $KMS_\infty$-states of $\mathcal A _K$ (see section \ref{extrkms}),
the verification of property $(5')$ is an immediate consequence of our construction and property $(6)$ follows by using Shimura's reciprocity law and the observation made in Prop. \ref{mainprop} (see section \ref{proof} for the details). \\ \\
Our paper is organized along the lines of this section, recalling on the way the necessary background. In addition we put some effort in writing a long Appendix which covers hopefully enough information to make this paper "readable" for a person which has little beforehand knowledge of the arithmetic theory of Shimura varieties.

\subsection*{Acknowledgment}
The author would like to thank James Milne, Sergey Neshveyev and Frederic Paugam for helpful and valuable comments
and his advisor Eric Leichtnam for his guidance through this first project of the author's Phd thesis.  \\
This work has been supported by the Marie Curie Research Training Network in Noncommutative Geometry, EU-NCG.

\subsection*{Notations and conventions}
We use the commen notations $\mathbb N$, $\ganz$,  $\q$,  $\real$, $\complex$. If $A$ denotes a ring or monoid, we denote its group of multiplicative units by $A^\times$. A number field is a finite extension of $\q$. The ring of integers of a number field $K$ is denoted by $\mathcal O_K$.
We denote by $\mathbb A_K = \mathbb A _{K,f} \times \mathbb A _{K,\infty}$ the adele ring of $K$ (with its usual topology), where $\mathbb A_{K,f}$ denotes the finite adeles and $\mathbb A _{K,\infty}$ the infinite adeles of $K$. $\mathbb A _K$ contains $K$ by the usual diagonal embedding and by $\widehat {\mathcal O}_K$ we denote the closure of $\mathcal O _K$ in $\mathbb A _{K,f}$. Invertible adeles are called ideles. The idele class group $\mathbb A_K^\times / K^\times$ of $K$ is denoted by $C_K$, its connected component of the identity by $D_K$. \\
We fix an algebraic closure $\qbar$ of $\q$ in $\complex$. Usually we think of a number field $K$ as lying in $\complex$ by an embedding $\tau : K \to \qbar \subset \complex$. Complex conjugation on $\complex$ is denoted by $\iota$. Sometimes we write $z^\iota$ for the complex conjugate of a complex number $z$.\\
Artin's reciprocity map $\mathbb A_K^\times \to Gal(K^{ab}/K) : \nu \mapsto [\nu]$ is normalized such that an uniformizing parameter maps to the arithmetic Frobenius element. Further given a group $G$ acting partially on a set $X$ we denote by \eq{\label{notlln}G \boxplus X = \{ (g,x) \in G\times X \ | \ gx \in X \}} the corresponding groupoid (see p. 327 \cite{Sergey}). \\
If X denotes a topological space we write $\pi_0(X)$ for its set of connected components.

%%%%%%%%%%%%%%%%%%%%%%%%%%%%%%%%%%%%%%%%%%%%%
%%%%%%%%        NUMBERTHEORETIC PART
%%%%%%%%%%%%%%%%%%%%%%%%%%%%%%%%%%%%%%%%%%%%%

\section{On the arithmetic subalgebra for $K=\q(i)$}
\label{bigexample}
Before we describe our general construction we will explain the easiest case $K=\q(i)$, where many simplifications occur, in some detail and point out the modifications necessary for the general case. For the reminder of this section $K$ always denotes $\q(i)$, although many of the definitions work in general. For the convenience of the reader we will try to make the following section as self-contained as possible. 

\subsection{The quotient map $U_{K} \to Z_{K}$}

We denote by $T^K$ the $\q$-algebraic torus given by the Weil restriction $T^K=Res_{K/ \q}(\mathbb G _{m,K})$ of the multiplicative group $\mathbb G_{m,K}$, i.e. for a $\q$-algebra $R$ the $R$-points of $T^K$ are given by $T^K(R) = (R \otimes_\q K)^\times$. In particular we see that $T^K(\q)=K^\times$, $T^K(\adfin) = \mathbb A_{K,f}^\times$ and $T^K(\real) = \mathbb A _{K,\infty}^\times$. In our special case we obtain that after extending scalars to $\real$ the $\real$-algebraic group $T^K_\real$ is isomorphic to $\mathbb S = Res_{\complex / \real}(\mathbb G_{m,\real})$. Further the finite set $X_K = T^K(\real) / T^K(\real)^+ = \pi_0(T^K(\real))$ consists in our case of only one point. With this in mind we consider the $0$-dimensional Shimura datum (see \ref{0shvD}) \eq{\mathcal S _K = (T^K,X_K,h_K)}
where the morphism $h_K : \mathbb S \to T^K_\real$ is simply given by the identity (thanks to $T^K_\real \cong \mathbb S $). (In the general case $h_K$ is chosen accordingly to Lemma \ref{lemmahodge}.)\\
The ($0$-dimensional) Shimura variety $Sh(\mathcal S_K)$ is in our case of the simple form \eq{Sh(\mathcal S _K)=T^K(\q)\backslash (X_K \times T^K(\adfin)) = K^\times \backslash \mathbb A_{K,f}^\times.}
We write $[z,l]$ for an element in $Sh(\mathcal S_K)$ meaning that $z \in X_K$ and $l \in T^K(\adfin)$. \\ 
(For general number fields the description of $Sh(\mathcal S _K)$ is less explicit but no difficulty occurs.)  
\begin{remark}
\label{remexgal}
The reader should notice that by class field theory we can identify $Sh(\mathcal S _K) = K^\times \backslash \mathbb A_{K,f}^\times$ with the Galois group $Gal(K^{ab}/K)$ of the maximal abelian extension $K^{ab}$ of $K$. This is true in general, see section \ref{symmbcmK}.
\end{remark}
The (topological) groupoid $U_K$ underlying the BC-type system $\mathcal A _K = (A_K,(\sigma_t)_{t\in \real})$ is now of the form (see (\ref{notlln}) for the notation) \eq{U_K = T^K(\adfin) \boxplus (\widehat{\mathcal O}_K \times Sh(\mathcal S _K) ) } with the natural action of $T^K(\adfin)$ on $Sh(\mathcal S _K)$ (see \ref{shvD}) and the partial action of $T^K(\adfin) = \mathbb A _{K,f}^\times$ on the multiplicative semigroup $\widehat{\mathcal O}_K \subset \mathbb A _{K,f}$ by multiplication. The group \eq{\Gamma_K^2 = \widehat {\mathcal O}_K^\times \times \widehat {\mathcal O}_K^\times} is acting on $U_K$ as follows
\eq{\label{standardaction}(\gamma_1,\gamma_2) (g,\rho,[z,l]) = (\gamma_1^{-1}g\gamma_2,\gamma_2 \rho,[z,l\gamma_2^{-1}]),}
where $\gamma_1,\gamma_2 \in \widehat{\mathcal O }_K^\times$, $g,l \in T^K(\adfin)$, $\rho \in \widehat{\mathcal O}_K$ and $z\in X_K$ and we obtain the quotient map \eq{U_K \longrightarrow Z_K = \Gamma_K^2 \backslash U_K.}
In the end of this section we will construct the arithmetic subalgebra $A_K^{arith}$ of the BC-type system $\mathcal A _K$ which is contained in $H_K=C_c(Z_K) \subset A_K$. For this we will need 

\subsection{The quotient map $U_{\siegel}\longrightarrow Z_\siegel$}

 In our case of $K=\q(i)$ the maximal CM subfield $E$ of $K$ is equal to $K$. The Shimura datum $\mathcal S_\siegel$ associated with $E$ is of the form (see \ref{construction}) \eq{\mathcal S _\siegel = (GSp(V_E,\psi_E),\mathbb H^±,h_{cm}).}
 Here $GSp = GSp(V_E,\psi_E)$ is the general symplectic group (cf. \ref{exD}) associated with the symplectic vector space $(V_E,\psi_E)$. \\
The latter is in general chosen accordingly to (\ref{injmap}). Due to the fact that the reflex field $E^*$ (cf. \ref{reflexB})
is equal to $E$ and the Serre group $S^E$ is equal to $T^E=T^K$ we can simply choose the $\q$-vector space $V_{E}$ to be the $\q$-vector space $E$ and the symplectic form $\psi_E : E\times E \to \q$ to be the map $(x,y) \mapsto Tr_{E/\q}(ixy^\iota)$. A simple calculation shows that $\psi_E(f(x),f(y))=det(f)\psi_E(x,y)$, for all $f \in End_\q(V_E)$ and all $x,y \in V_E$, therefore we can identify $GSp$ with $GL_2=GL(V_E)$. Now again using the fact that the Serre group $S^E$ equals $T^E$ we see that the general construction of $h_{cm} : \mathbb S = T^E_\real \to GSp_\real = GL_{2,\real}$ (see (\ref{defhcm})) is given on the $\real$-points by $a+ib \in \complex^\times =  \mathbb S(\real) \to \matrt{a& -b \\ b & a} \in GL_2(\real)$. Each $\alpha \in GSp(\real)$ defines a map $\alpha^{-1}h_{cm}\alpha : \mathbb S  \to GSp_\real$ given on the $\real$-points by $a+ib \in \complex^\times \mapsto \alpha^{-1}h_{cm}(a+ib)\alpha \in GL_2(\real)$ and the $GSp(\real)$-conjugacy class $X = \{\alpha^{-1}h_{cm}\alpha \  | \  \alpha \in GSp_\real(\real)\}$ of $h_{cm}$ can be identified with the Siegel upper lower half space $\mathbb H^± = \complex - \real$ by the map \eqst{\alpha^{-1}h_{cm}\alpha \in X \mapsto (\alpha^{-1}h_{cm}(i)\alpha) \cdot i \in \mathbb H^±,} where the latter action $\cdot$ denotes Moebius transformation. Under this identification the morphism $h_{cm}$ corresponds to the point $x_{cm}=i$ on the upper half plane $\mathbb H$. The point $x_{cm} \in \mathbb H$ is a so called CM-point (see \ref{canmodD}).  
\begin{remark}
The definition of a CM-point and the observation made in (\ref{emb}) explain the need of using the Serre group in the general construction of $h_{cm}$. The explanation given in \ref{intro2} shows in particular why we have to define the vector space $V_E$ in general accordingly to (\ref{injmap}). 
\end{remark}
The Shimura variety $Sh(\mathcal S_\siegel)$ is of the nice form (cf. (\ref{sv5})) \eqst{Sh(\mathcal S _\siegel) = GSp(\q)\backslash (\mathbb H^± \times GSp(\adfin) ).}
Again we write elements as $[z,l] \in Sh(\mathcal S_\siegel)$ with $z \in \mathbb H^±$ and $g \in GSp(\adfin)$.\\
In our case the topological groupoid $U_\siegel$ underlying the Shimura system $\mathcal A _\siegel$ is given by (cf. \ref{costumeII})
\eq{U_\siegel = GSp(\adfin) \boxplus (M_2(\widehat{\ganz}) \times Sh(\mathcal S _\siegel )),}
where $GSp(\adfin) = GL_2(\adfin)$ is acting in the natural way on $Sh(\mathcal S_\siegel)$ and partially on the multiplicative monoid of $2\times2$-matrices $M_2(\widehat \ganz) \subset M_2(\adfin)$ with entries in $\widehat \ganz$. The group \eq{\Gamma_\siegel^2 = GL_2(\widehat \ganz) \times GL_2(\widehat \ganz)} is acting on $U_\siegel$ exactly like in ($\ref{standardaction}$) and induces the quotient map \eq{U_\siegel \longrightarrow Z_\siegel = \Gamma_\siegel^2 \backslash U_\siegel.}
\begin{remark}
Note that the quotient $Z_\siegel$ is not a groupoid anymore (see top of p. 251 \cite{HaPa}).
\end{remark}
In our example it is sufficient to consider the positive groupoid $U_\siegel^+$ (see \ref{somequotientmaps}) associated with $U_\siegel$. (In the general case the adjoint groupoid $U_\siegel^{ad}$ seems to be more appropriate.) It is given by \eq{U_\siegel^+ = GSp(\q)^+ \boxplus (M_2(\widehat \ganz) \times \mathbb H)}
together with the group \eq{(\Gamma_\siegel^+)^2 = GL_2(\ganz)^+ \times GL_2(\ganz)^+ = SL_2(\ganz)\times SL_2(\ganz)} acting by \eq{(\gamma_1,\gamma_2)(g,\rho,z)=(\gamma_1 g \gamma_2^{-1},\gamma_2 \rho,\gamma_2 z)} and inducing the quotient map \eq{\label{exZplus}U^+_\siegel \longrightarrow Z_\siegel^+ = (\Gamma_\siegel^+)^2 \backslash U^+_\siegel.} By construction $GSp(\real)$ is acting (free and transitively) on $\mathbb H^±$ and $GSp(\real)^+$, the connected component of the identity, can be thought of as stabilizer of the upper half plane $\mathbb H ^+ = \mathbb H$ which explains the action of $GSp(\q)^+ = GSp(\q) \cap GSp(\real)^+$ on $\mathbb H$. Thanks to criterion \ref{crit} we know that the natural (equivariant) morphism of topological groupoids $U^+_\siegel \to U_\siegel$ given by $(g,\rho,z) \mapsto (g,\rho,[z,1])$ induces a homeomorphism on the quotient spaces $Z_\siegel^+ \longrightarrow Z_\siegel$ in the commutative diagram
\diagramlab{critex1}{U_\siegel \ar[rr] && Z_\siegel \\ U_\siegel^+ \ar[rr] \ar[u] && Z_\siegel^+ \ar[u]_{\cong}.}

\begin{remark}
The groupoid $U_\siegel^+$ corresponds to the $GL_2$-system of Connes and Marcolli (see \cite{conmar} and 5.8 \cite{HaPa})
\end{remark}

\subsection{A map relating $U_K$ and $U_\siegel$}

We want to define an equivariant morphism of topological groupoids $U_K \longrightarrow U_\siegel$, where equivariance is meant with respect to the actions of $\Gamma_K$ on $U_K$ and $\Gamma_\siegel$ on $U_\siegel$.   
For this it is necessary (and more or less sufficient) to construct a morphism of Shimura data between $\mathcal S_K = (T^K,X_K,h_K)$ and $\mathcal S _\siegel = (GSp(V_E,\psi_E)),\mathbb H^±,h_{cm})$, which is given by a morphism of algebraic groups \eq{\varphi: T^K \longrightarrow GSp} such that $h_{cm} = \varphi_\real \circ h_{K}$. In our case the general construction of $\varphi$, stated in (\ref{defphi}), reduces to the simple map (on the $\q$-points) \eqst{a+ib \in K^\times=T^K(\q) \mapsto \matrt{a & -b \\ b & a} \in GL_2(\q)=GSp(\q)}
and we see in fact that after extending scalars to $\real$ the morphism $\varphi_\real : T^K_\real = \mathbb S \to GSp_\real$ is already equal to $h_{cm}:\mathbb S \to GSp_\real$. The simplicity of our example comes again from the fact that we don't have to bother about the Serre group, which makes things less explicit, although the map $\varphi : T^K \to GSp$ still has a quite explicit description even in the general case thanks to Lemma \ref{lemmaexplicitphi}. \\ \\
Now by functoriality (see section \ref{0shvD}) we obtain a morphism of Shimura varieties \eqst{Sh(\varphi):Sh(\mathcal S_K ) \to Sh(\mathcal S _\siegel)} which can be explicitly described by \alst{[z,l] \in K^\times \backslash (X_K \times T^K(\adfin)) \mapsto [x_{cm},\varphi(\adfin)(l)] \in GSp(\q)\backslash(\mathbb H^± \times GSp(\adfin)).}
In the general case we have essentially the same description (see (\ref{mapbcmgroupoid})), the point being that every element $z$ in $X_K$ is mapped to $x_{cm} \in \mathbb H^±$, as in the general case. Using $\widehat{\mathcal O}_K = \widehat \ganz \otimes _\ganz \mathcal O_K$ we can continue the map $\varphi(\adfin)$ to a morphism of (topological) semigroups $M(\varphi)(\adfin): \widehat {\mathcal O}_K \to M_2(\widehat \ganz)$ by setting $n \otimes (a+ib) \mapsto \matrt{an & -bn \\ bn & an}$. By continuation we mean that $\varphi(\adfin)$ and $M(\varphi)(\adfin)$ agree on the intersection of $T^K(\adfin) \cap \widehat{\mathcal O}_K \subset \mathbb A_{K,f}$. In the general case the explicit description of $\varphi$ given in (\ref{phiexplicit}) is used to continue $\varphi$ to $M(\varphi)$ (see \ref{continuation}), the above example being a special case. Now it can easily be checked that \eq{(g,\rho,[z,l])\in U_K \mapsto (\varphi(\adfin)(g),M(\varphi)(\adfin)(\rho),[x_{cm},\varphi(\adfin)(l)]) \in U_\siegel}
defines the desired equivariant morphism of topological groupoids \eq{U_K \longrightarrow U_\siegel.}
Summarizing we obtain the following commutative diagram (using (\ref{critex1}))
\diagram{U_K \ar[rr] \ar[d] && Z_K \ar[d] \\ U_\siegel \ar[rr] && Z_\siegel \\ U^+_\siegel \ar[u] \ar[rr] && Z^+_\siegel \ar[u]_{\cong}}
which gives us the desired morphism of topological spaces \eq{\label{thetaex}\Theta : Z_K \longrightarrow Z_\siegel \longrightarrow Z^+_\siegel.}
In general we have to go one step further and use the adjoint groupoid $Z^{ad}_\siegel$ (cf. \ref{somequotientmaps}) but this only due to the general description of the automorphism group $Aut_\q(\mathcal M)$ of the field of arithmetic automorphic functions $\mathcal M$, see \ref{defmathcalM} and the next section for explanations. 

\subsection{Interlude: Theory of Complex Multiplication}
\label{exgenab}
In this section we will provide the number theoretic background which is necessary to understand the constructions done so far. \\
We are interested in constructing the maximal abelian extension $E^{ab}$ of $E = K = \q(i)$ and there are in general two known approaches to this problem. 

\subsubsection*{I) The elliptic curve $A: y^2=x^3+x$}

Let us denote by $A$ the elliptic curve defined by the equation \eq{\label{ellcur}A : y^2=x^3+x.} It is known that the field of definition of $A$ and its torsion points generate the maximal abelian extension $E^{ab}$ of $E$. (By field of definition of the torsion points of $A$ we mean the coordinates of the torsion points.) \\
Notice that in our case the complex points of our elliptic curve $A$ are given by $A(\complex) = \complex / \mathcal O_E$ and the rational ring of endomorphisms of $A$ turns out to be $End(A(\complex))_\q = \q\otimes_\ganz \mathcal O_E = E$,we say that $A$ has complex multiplication by $E$. Remember that $\mathcal O _E = \ganz[i]$.
\begin{remark}
To obtain the abelian extensions of $K$ provided by $A$ explicitly one may use for example the Weierstrass $\mathfrak p$-function associated with $A$ (see \cite{silver}) but as we will use another approach we don't want to dive into this beautiful part of explicit class field theory.
\end{remark}
\subsubsection*{II) The (Siegel) Modular curve $Sh(GSp,\mathbb H^±,h_{cm})$}

We want to interpret the Shimura variety $Sh(\mathcal S _\siegel) = Sh(GSp,\mathbb H^±,h_{cm})$ constructed in the last section as moduli space of elliptic curves with torsion data. 

\subsubsection*{The moduli theoretic picture}

For this we consider the connected component $Sh^o = Sh(\mathcal S _\siegel)^o$ of our Shimura variety which is described by the projective system (see \ref{cshvD} and \cite{isv}) \eq{Sh^o = \varprojlim_N \Gamma(N) \backslash \mathbb H,}
where $\Gamma(N)$, for $N\geq 1$, denotes the subgroup of $\Gamma = \Gamma(1) = SL_2(\ganz)$ defined by $\Gamma(N) = \{ g \in \Gamma \ | \ g \equiv \matrt{1&0\\0&1} \text{ mod } N \}$. We can view the quotient $\Gamma(N) \backslash \mathbb H$ as a complex analytic space, but thanks to the work of Baily and Borel, it carries also a unique structure of an algebraic variety over $\complex$ (see \cite{isv}). We will use both viewpoints. Seen as an analytic space we write $\mathbb H(N) = \Gamma(N) \backslash \mathbb H$ and for the algebraic space we write $Sh^o_N = \Gamma(N) \backslash \mathbb H$.  \\ \\
Observe now that the space $\mathbb H(N)$ classifies isomorphism classes of pairs $(A,t)$ given by an elliptic curve $A$ over $\complex$ together with a $N$-torsion point $t$ of $A$. In particular $\mathbb H (1) = \Gamma \backslash \mathbb H$ classifies isomorphism classes of elliptic curves over $\complex$. \\ \\
In this picture our CM-point $[x_{cm}]_1=[i]_1 \in \mathbb H (1)$ corresponds to the isomorphism class of the elliptic curve $A$ from (\ref{ellcur}) and more general the points $[x_{cm}]_N \in \mathbb H(N)$ capture the field of definition of $A$ and its various torsion points and recover the maximal abelian extension $E^{ab}$ of $E$ in this way! By $[z]_N$ we denote the image of $z \in \mathbb H$ in $\mathbb H(N)$ under the natural quotient map $\mathbb H \to \mathbb H(N)$.

\begin{remark}
For the relation between $Sh$ and $Sh^o$ we refer the reader to pp. 51 \cite{isv}.
\end{remark}

\subsubsection*{The field of arithmetic Modular functions $\mathcal M$}

To construct the abelian extensions provided by the various points $[x_{cm}]_N$ explicitly we proceed as follows: We consider the connected canonical model $M^o$ of $Sh^o$ (see \ref{ccanmodD}) which provides us with an algebraic model $M^o_N = \Gamma(N) \backslash M^o$ of the algebraic variety $Sh^o_N$ over the cyclotomic field $\q(\zeta_N)$. In general we obtain algebraic models over subfields of $\q^{ab}$. This means $M^o_N$ is an algebraic variety defined over the cyclotomic field $\q(\zeta_N)$ and after scalar extension to $\complex$ it becomes isomorphic to the complex algebraic variety $Sh^o_N$. Let us denote by $k(M^o_N)$ the field of rational functions on $M^o_N$, in particular this means elements in $k(M^o_N)$ are rational over $\q(\zeta_N)$. It makes sense to view the point $[x_{cm}]_N$ as a point on $M^o_N$ and if a function $f \in k(M^o_N)$ is defined at $[x_{cm}]_N$, then we know (cf. \ref{defmathcalM}) \eq{ \label{abextex} f([x_{cm}]_N) \in E^{ab}.}
In particular varying over the various $N$ and the rational functions in $k(M^o_N)$ the values $f([x_{cm}]_N)$
generate $E^{ab}$ over $E$. The next step is to realize that the function field $k(M^o_N)$ can be seen as subset of the field of rational functions $k(Sh^o_N)$ on the complex algebraic variety $Sh^o_N$ (cf. \ref{defmathcalM}). As rational functions in $k(Sh^o_N)$ correspond to meromorphic functions on $\mathbb H (N)$ and meromorphic functions on $\mathbb H (N)$ are nothing else than meromorphic functions on $\mathbb H$ that are invariant under the action of $\Gamma(N)$ we can view each rational function in $k(M^o_N)$ as a meromorphic function on $\mathbb H$ which is invariant under $\Gamma(N)$. If we denote by $k(M^o_N)_{cusp}$ the subfield of $k(M^o_N)$ consisting of functions $f \in k(M^o_N)$ that give rise to meromorphic functions on $\mathbb H$ that are meromorphic at the cusps (see \ref{defmathcalM}), it makes sense to define the field of meromorphic functions $\mathcal M$ on $\mathbb H$ given by the union \eq{\mathcal M = \bigcup_N \  k(M^o_N)_{cusp}.}
Due to (\ref{abextex}) we know furthermore that for every $f \in \mathcal M$, which is defined in $x_{cm}$, we have \eq{f(x_{cm}) \in K^{ab} = E^{ab}} and $K^{ab}$ is generated in this way. Therefore we call the field $\mathcal M$ the \textbf{field of arithmetic Modular functions}. \\
Explicitly $\mathcal M$ is described for example in \cite{Shimura} or \cite{conmar} (Def. 3.60). A very famous arithmetic Modular function is given by the $j$-function which generates the Hilbert class field of an arbitrary imaginary quadratic field.
\begin{remark}
1) In light of \ref{intro2} and \ref{proph} 2) we mention that the field of definition $E(x_{cm})$ of $x_{cm}$ is in our example equal to $E= \q(i)$, which is the reason why our example is especially simple.\\
2) If we take a generic meromorphic function $g$ on $\mathbb H (N)$ that is defined in $[x_{cm}]_N$ then the value $g([x_{cm}]_N) \in \complex$ will not even be algebraic. 
This is the reason why we need the canonical model $M^o$ which provides an artithmetic structure for the field of meromorphic functions on $\mathbb H(N)$.
\end{remark}
In the general case the construction of $\mathcal M$ is quite similar to the construction above (see section \ref{cmtheory}), the only main difference being that in general the theory is much less explicit (e.g. the description of $\mathcal M$).

\subsubsection*{Automorphisms of $\mathcal M$ and Shimura's reciprocity law}

In our example $K=E=\q(i)$, using the notation from \ref{goingdown}, we have the equality $\overline{\mathcal E} = \frac {GSp(\adfin)} {\q^\times}$ and obtain a group homomorphism \eq{GSp(\adfin) \longrightarrow \overline {\mathcal E} \longrightarrow  Aut_\q(\mathcal M),}
where the first arrow is simply the projection and the second arrow comes from \ref{autmathcalM}. We denote the action of $\alpha \in GSp(\adfin)$ on a function $f \in \mathcal M$ by ${}^{\alpha}\hspace{-0.6mm}f$. In particular we see that $\alpha \in GSp(\q)^+ = SL_2(\q)$ is acting by \eq{\label{loione}{}^{\alpha}\hspace{-0.6mm}f = f \circ \alpha^{-1},}
where $\alpha^{-1}$ acts on $\mathbb H$ by Mobius transformation. The adjoint system, which is well suited for the general case (cf. \ref{somequotientmaps}), is not needed in our special example. We note that the group $GSp^{ad}(\q)^+$ occuring in (\ref{somequotientmaps}) is given by $SL_2(\q) / \{±I\}$, where $I$ denotes the unit in $SL_2(\q)$. Due to the fact that $\{±I\}$ acts trivially on $\mathbb H$ we can lift the action to $GSp(\q)^+$.
Further there is a morphism of algebraic groups \eq{\eta : T^K \to GSp} 
which induces a group homomorphism denoted by (compare (\ref{overlineeta})) \eq{\overline{\eta} = \eta(\adfin): T^K(\adfin) \to GSp(\adfin).} The \textbf{reciprocity law of Shimura} can be stated in our special case as follows: \\
Let $\nu$ be in $\mathbb A_{K,f}^\times$, denote by $[\nu] \in Gal(K^{ab}/K)$ its image under Artin's reciprocity map and let $f \in \mathcal M$ be defined in $x_{cm} \in \mathbb H$. Then ${}^{\overline{\eta}(\nu)}\hspace{-0.6mm}f $ is also defined in $x_{cm} \in \mathbb H$ and \al{\label{loideux}{}^{\overline{\eta}(\nu)}\hspace{-0.6mm}f(x_{cm}) = [\nu]^{-1}f(x_{cm}) \in K^{ab}.}
The formulation in the general case concentrates on the group $\overline{\mathcal E}$ (see \ref{recilawgeneral}). \\ \\
Now we are ready for the

\subsection{Construction of the arithmetic subalgebra}

Define $\mathcal M^{cm}$ to be the subring of functions in $\mathcal M$ which are defined at $x_{cm} \in \mathbb H$.
For every $f \in \mathcal M^{cm}$ we define a function $\widetilde f $ on the groupoid $U_\siegel^+$ by \al{\widetilde f (g,\rho,z) =   \left\{ \begin{array}{cc} ^{\rho}\hspace{-0.5mm}f(z)& \rho \in GSp(\widehat \ganz)  \\ 0&  \rho \in M_2(\widehat \ganz) - GSp(\widehat \ganz) \end{array} \right.}
Thanks to (\ref{loione}) $\widetilde f$ is invariant under the action of $\Gamma_\siegel^+ = SL_2(\ganz)\times SL_2(\ganz)$ and therefore $\widetilde f$ descends to the quotient $Z_\siegel^+$ (cf. (\ref{exZplus})). Proposition \ref{pullbackprop} shows further that the pull-back $\widetilde f \circ \Theta$ (see (\ref{thetaex})) defines a compactly supported, continuous function on $Z_K$, i.e. we have $\widetilde f \circ \Theta \in H_K = C_c(Z_K) \subset A_K$. Therefore we can define the $K$-subalgebra $A^{arith}_K$ of $H_K$ generated by these elements
\eq{A^{arith}_K = <\widetilde f \circ \Theta \ | \ f \in \mathcal M ^{cm}>_K.} 
Now we want to show that $A^{arith}_K$ is indeed an arithmetic subalgebra for $\mathcal A _K$ (see (\ref{arithalgdef})). \\
The set $\mathcal E _\infty$ of extremal $KMS_\infty$-states of $\mathcal A _K$ is indexed by the set $Sh(\mathcal S_K)$ and for every $\omega \in Sh(\mathcal S_K)$ the corresponding $KMS_\infty$-state $\varrho_\omega$ is given on an element $f \in H_K = C_c(Z_K)$ by (cf. \ref{extrkms}) \eq{\varrho_\omega(f) = f(1,1,\omega).}
Using remark \ref{remexgal} we write $[\omega]$ for an element $\omega \in Sh(\mathcal S _K)$ when regarded as element in $Gal(K^{ab}/ K)$.
Now if we take a function $f \in \mathcal M ^{cm}$ and $\omega \in Sh(\mathcal S_K)$ we immediately see \eq{\varrho_\omega(\widetilde f \circ \Theta) = [\omega]^{-1}(f(x_{cm})) \in K^{ab}}
and property (v) follows (in general we will only show property (v') of course). To show property (vi) we take a symmetry $\nu \in C_K = \mathbb A_K^{\times} / K^\times$ (see \ref{symmbcmK}) of $\mathcal  A_K$ and denote by $[\nu] \in Gal(K^{ab}/K)$ its image under Artin's reciprocity homomorphism and let $f$ and $\omega$ be as above. We denote the action (pull-back) of $\nu$ on $\varrho_\omega$ by $^{\nu}\hspace{-0.5mm}\varrho_\omega$ and obtain \eq{^{\nu}\hspace{-0.5mm}\varrho_\omega(\widetilde f \circ \Theta) = {}^{\varphi(\adfin)(\nu)}\hspace{-0.5mm}f(x_{cm}).}
But thanks to our observation Proposition \ref{mainprop} we know that \eq{\varphi(\adfin)(\nu) = \overline\eta (\nu) \in GSp(\adfin)} and now thanks to Shimura's (\ref{loideux}) we can conclude \eq{^{\nu}\hspace{-0.5mm}\varrho_\omega(\widetilde f \circ \Theta) = {}^{\overline \eta (\nu)}\hspace{-0.5mm}f(x_{cm}) = [\nu]^{-1}(f(x_{cm}))=[\nu]^{-1}(\varrho_\omega(\widetilde f \circ \Theta)) \in K^{ab}}
which proves property (vi). For the general case and more details we refer the reader to \ref{proof}. 

\begin{remark}
Our arithmetic subalgebra $\mathcal A _K$ in the case of $K = \q(i)$ is essentially the same as in \cite{cmr1}.
\end{remark}
\begin{remark}
In a fancy (and very sketchy) way we might say that the two different pictures, one concentrating on the single elliptic curve $A$ the other on the moduli space of elliptic curves (see the beginning of section \ref{exgenab}), are related via the Langlands correspondence. In terms of Langlands correspondence, the single elliptic curve $A$ lives on the motivic side whereas the moduli space of elliptic curves lives (partly) on the automorphic side. As we used the second picture for our construction of an arithmetic subalgebra, we might say that our construction is automorphic in nature. This explains the fact that we have a "natural" action of the idele class group on our arithmetic subalgebra (see above). Using the recent theory of endomotives (see chapter 4 and in particular p. 551 \cite{conmar}) one can recover the arithmetic subalgebra $A_K^{arth}$ by only using the single elliptic curve $A$. In particular one obtains a natural action of the Galois group. This and more will be elaborated in another paper.
\end{remark}
We now concentrate on the general case.

\section{Two Shimura data and a map}
\label{shimdata}
As throughout the paper let $K$ denote a number field and $E$ its maximal CM subfield. We fix an embedding $\tau : K \to \qbar \to \complex$ and denote complex conjugation on $\complex$ by $\iota$. \\
To $K$, resp. $E$, we will attach a Shimura datum $\mathcal S_K$, resp. $\mathcal S_{\siegel}$, and show how to construct a morphism $\varphi : \mathcal S_K \to \mathcal S_{\siegel}$ between them. We will freely use the Appendix: every object not defined in the following can be found there or in the references given therein. \\Recall that the Serre group attached to $K$ is denoted by $S^K$ (cf. \ref{defC}), it is a quotient of the algebraic torus $T^K$ (defined below), the corresponding quotient map is denoted by $\pi^K : T^K\to S^K$.  

\subsection{Protagonist I: $\mathcal S_K$}
\label{shdatum1}
The $0$-dimensional Shimura datum $\mathcal S_K=(T^K,X_K, h_K)$ (see section \ref{0shvD}) is given by the Weil restriction $T^K=Res_{K/\q}(\mathbb G_{m,K})$, the discrete space $X_K = T^K(\real) / T^K(\real)^+$ and a morphism $h_K:\mathbb S = Res_{\komplex / \real}(\mathbb G_{m,\complex})\to T^K_\real$ which is chosen accordingly to the next
\begin{lemma}
\label{lemmahodge}
There is a morphism of algebraic groups $h_K : \mathbbm S \to T^K_\reell$ such that the diagram
\diagram{\mathbbm S \ar[r]^{h_K} \ar[dr]^{h^K} & T^K_\reell \ar[d]^{\pi^K_\real} \\ & S^K_\real }
commutes.
\end{lemma}
\begin{proof}
Remember that $h^K: \mathbb S \to S^K_\real$ is defined  as the composition \diagram{\mathbb S \ar[rrr]^{Res_{\komplex / \reell}(\mu^K) \ \ \ \ \ \ \ \ } & & & Res_{\komplex / \reell}(S^K_\komplex) \ar[rr]^{\ \ \ \ \ \ \ Nm_{\komplex / \reell}}  & & S^K_\reell,}
where $\mu^K : \mathbb G _{m,\complex} \to T^K_\complex$ is defined by $\mu^K = \pi^K_\complex \circ \mu_\tau$ (cf. \ref{canmapC}).
Define $h_K : \mathbb S \to T^E_\reell$ simply by 
\diagram{\mathbb S \ar[rrr]^{Res_{\komplex / \reell}(\mu_\tau) \ \ \ \ \ \ \ \ } & & & Res_{\komplex / \reell}(T^K_\komplex) \ar[rr]^{\ \ \ \ \ \ \ Nm_{\komplex / \reell}}  & & T^K_\reell.}
For proving our claim it is enough to show that the following diagram
\diagram{
\mathbb S \ar[rrrrd]_{Res_{\komplex / \reell}(\mu^K)\ \ \ \ \ \ } \ar[rrrr]^{Res_{\komplex / \reell}(\mu_\tau) \ \ \ \ } & & & & Res_{\komplex/ \real}(T^K_\komplex) \ar[rr]^{\ \ \ \ \ \  Nm_{\komplex / \reell}} \ar[d]^{Res_{\komplex / \real}(\pi_\komplex^K)} & & T^K_\reell \ar[d]^{\pi_\real^K}
\\
 & & & & Res_{\komplex/ \real}(S^K_\komplex) \ar[rr]^{\ \ \ \ \ \ Nm_{\komplex / \reell}} & & S^K_\reell
}
is everywhere commutative. \\ The triangle on the left is commutative, because $Res_{\complex / \real}$ is a functor. Thanks to theorem \ref{thmA} it is enough to show that rectangle on the right is commutative after applying the functor $X^*$ (cf. \ref{charA}). Since $\pi^K : T^K \to S^K$ is defined to be the inclusion $X^*(S^K) \subset X^*(T^K)$ on the level of characters (cf. \ref{defC}), we see that $X^*(Res_{\complex / \real} (\pi^K_\complex) )$ and $X^*(\pi^K_\real)$ are inclusions as well and the commutativity follows.

\end{proof}
 
\subsection{Protagonist II: $\mathcal S_{\siegel}$}
\label{construction}
The construction of the Shimura datum $\mathcal S_{\siegel}$ in this section goes back to Shimura \cite{Shimura}, see also \cite{wei}. 
It is of the form $\mathcal S_{\siegel} = (GSp(V_E,\psi_E),\mathbb H^±_g,h_{cm})$. 
The symplectic $\q$-vector space $(V_E,\psi_E)$ is defined as follows. \\
Choose a finite collection of primitive CM types $(E_i,\Phi_i)$, $1\leq i \leq r$, such that \\
i) for all $i$ the reflex field $E_i^*$ is contained in $E$, i.e. $\forall i \ E_i^* \subset E$, and \\
ii) the natural map (take (\ref{muphiB}) and apply the universal property from \ref{defC})
\begin{equation}
\label{injmap}
\xymatrix{S^E \ar[rr]^{ \prod N_{E/E_i^*} \ \ \  } & & \prod_{i=1}^r S^{E_i^*} \ar[rr]^{\prod \rho_{\Phi_i}} & & \prod_{i=1}^r T^{E_i}}
\end{equation} is injective. (Proposition 1.5.1 \cite{wei} shows that this is allways possible.) \\
For every $i \in \{1,..,r\}$ we define a symplectic form $\psi_i : E_i \times E_i \to \q$ on $E_i$ by choosing a totally imaginary generator $\xi_i$ of $E_i$ (over $\q$) and setting \eq{\label{defpsi}\psi_i(x,y) = Tr_{E_i / \q}(\xi_i x  y^\iota).}
Now we define $(V_E,\psi_E)$ as the direct sum of the symplectic spaces $(E_i,\psi_i)$. Instead of $GSp(V_E,\psi_E)$ we will sometimes simply write $GSp$. \\ \\
To define the morphism $h_{cm}$ the essential step is to observe (see Remark 9.2 \cite{pedestrian}) that the image of the map $\rho_{\Phi_i} \circ N_{E/E_i^*} : S^E \to T^{E_i}$ is contained in the subtorus $\mathcal T^{E_i}$ of $T^{E_i}$, which is defined on the level of $\q$-points by 
\eq{\mathcal T ^{E_i}(\q) = \{ x \in E_i^\times \ | \ x x^\iota \in \q^\times \}}
and analogously $T^{E_i}(R)$ is defined for an arbitrary $\q$-algebra $R$. This is an important observation, because there is an obvious inclusion of algebraic groups (cf. \ref{exA}) \eq{\label{emb}\mathfrak i :  \prod_{i=1}^r \ \mathcal T^{E_i} \to GSp(V_E,\psi_E),}
whereas there is in general \underline{\textbf{no}} embedding $\prod T^{E_i} \to GSp$.\\
With this in mind we define $h_{cm}$ as the composition \diagramlab{defhcm}{\mathbbm S \ar[r]^{h^E} & S^E_\real \ar[rr]^{\prod N_{E / E_i^*,\real} \ \ \ \ \ } & & \prod_{i=1}^r S^{E_i^*}_\real \ar[rr]^{\prod \rho_{\Phi_i,\real} } & & \prod_{i=1}^r  \mathcal T^{E_i}_\real \ar[r]^{\mathfrak i_\real } & GSp_\real.}
Write $h'_{cm} : \mathbb S \to \prod_{i=1}^r  \mathcal T^{E_i}_\real$ for the composition of the first three arrows.
\begin{remark}
1) By construction $h_{cm}$ is a CM point (cf. \ref{canmodD}) which is needed later to construct explicitly abelian extensions $K$. See \ref{intro2}.  \\
2) Viewed as a point on the complex analytic space $\mathbb H^±_g$ we write $x_{cm}$ instead of $h_{cm}$. Further we denote the connected component of $\mathbb H^±_g$ containing $x_{cm}$ by $\sieg $, i.e. $x_{cm} \in \sieg$.
\end{remark}
Our CM point $h_{cm}$ enjoys the following properties
\begin{lemma}
\label{proph}
1) We have $h_{cm} = \mathfrak i \circ \prod_{i=1}^r h_{\Phi_i}$ (see (\ref{hphiB})). \\
2) The field of definition $E(x_{cm})$ of $x_{cm}$ is equal to the composite of the reflex fields $\widetilde E = E_1^* \cdots E_r^* \subset E$, i.e. the associated cocharacter $\mu_{cm}$ of $h_{cm}$ is defined over $\widetilde E$ (see \ref{canmodD}).\\
3) The $GSp(\real)$-conjugacy classes of $h_{cm}$ can be identified with the Siegel upper-lower half plane $\mathbb H^±_g$, for some $g \in \mathbb N$ depending on $E$.
\end{lemma}
\begin{proof}
1) This follows immediately from (\ref{comphphi}) and (\ref{comphnorm}). \\
2) This follows from p. 105 \cite{isv} and 1). \\
3) For this we refer to the proof of lemma 3.11 \cite{wei}.
\end{proof}

\subsection{The map $\varphi : \mathcal S_K \to \mathcal S_{\siegel}$}
\label{themapphi}
On the level of algebraic groups $\varphi : T^K \to GSp$ is simply defined as the composition
\diagramlab{defphi}{T^K \ar[r]^{\pi^K} & S^K \ar[r]^{N_{K/E}} & S^E \ar[rr]^{\prod N_{E / E_i^*} \  \ \ \ } && \prod_{i=1}^r S^{E_i^*} \ar[r]^{\ \ \ \ \prod \rho_{\Phi_i} \ \ \ \ } &  \prod_{i=1}^r  \mathcal T^{E_i} \ar[r]^{\ \ \ \mathfrak i} & GSp.}
For $\varphi$ being a map between $\mathcal S_K$ and $\mathcal S_{\siegel}$ we have to check that the diagram
 \diagramlab{phihcm}{\mathbbm S \ar[r]^{h_K} \ar[dr]_{h_{cm}} & T^K_\real \ar[d]^{\varphi_\real} \\ & GSp_\real} commutes,
 but this compatibility is built into the construction of $h_K$.
 Using the reflex norm (cf. \ref{reflexB} and (\ref{reflexnormprop})) we can describe $\varphi$ as follows
 \begin{lemma}
 \label{lemmaexplicitphi}
 The map $\varphi : T^K \to GSp$ is equal to the composition \diagramlab{phiexplicit}{T^K \ar[rr]^{\prod N_{K/ E_i^*}} & & \prod T^{E_i^*} \ar[rr]^{\prod N_{\Phi_i}} & &  \prod_{i=1}^r  \mathcal T^{E_i} \ar[r]^{\ \ \ \mathfrak i} & GSp .}
 \end{lemma}

\section{About arithmetic Modular functions}
\label{cmtheory}
\subsection{Introduction}
\label{intro2}
We follow closely our references \cite{Del2}, \cite{MiS} and \cite{wei}. See also \cite{hida}. The reader should be aware of the fact that we are using a different normalization of Artin's reciprocity map than in \cite{MiS} and have to correct a "sign error" in \cite{Del2} as pointed out on p. 106 \cite{isv}.  \\ \\
As usual we denote by $K$ a number field containing a CM subfield and denote by $E$ the maximal CM subfield of $K$. In this section we want to explain how the theory of Complex Multiplication provides (explicit) abelian extensions of $K$. In general one looks at a CM-point $x \in X$ on a Shimura variety $Sh(G,X)$ and by the theory of canonical models one knows that the point $[x,1]$ on the canonical model $M(G,X)$ of $Sh(G,X)$ is rational over the maximal abelian extension $E(x)^{ab}$ of the field of definition $E(x)$ of $x$ (see \ref{canmodD}). \\ \\
In our case we look at Siegel modular varieties $Sh(GSp,\mathbb H^±_g)$ which can be considered as (fine) moduli spaces of Abelian varieties over $\complex$ with additional data (level structure, torsion data and polarization). See chap. 6 \cite{isv} for an explanation of this. Each point $x \in \mathbb H_g$ corresponds to an Abelian variety $A_x$. \\
In opposite to the case of imaginary quadratic fields, in general the field of definition $E(x)$ is neither contained in $E$ nor in $K$. Therefore, in order to construct abelian extensions of $K$, we have to find an Abelian variety $A_x$ such that \eq{E(x) \subset K.}
This is exactly the reason for our choice of $x_{cm}\in \mathbb H_g$ because we know (see \ref{proph} 2)) that 
\eq{E(x_{cm}) = \widetilde E = E_1^* \cdots E_r^* \subset E \subset K. } 
Here $x_{cm}$ corresponds to a product $A_{cm} = A_1 \times \cdots \times A_r$ of simple Abelian varieties $A_i$ with complex multiplication given by $E_i$. This construction is the best one can do to generate abelian extensions of $K$ using the theory of Complex Multiplication. The miracle here is again that the field of definition of $A_{cm}$ and of its torsion points generate abelian extensions of $E(x_{cm})$. Now to obtain these abelian extensions explicitly one proceeds in complete analogy with the case of $\q(i)$
explained in \ref{exgenab}, namely rational functions on the connected canonical model $M^o$ of the connected Shimura variety $Sh(GSp,\mathbb H^±_g)^o$ give rise to arithmetic Modular functions on $\mathbb H_g$ which generate the desired abelian extensions when evaluated at $x_{cm}$. This will be explained in detail in the following.

\subsection{Working over $\qbar$}

\subsubsection{The field $\mathcal F$ of arithmetic automorphic functions}

We start with the remark that the reflex field of $(GSp,\mathbb H^±_g)$ (cf. \ref{canmodD}) is equal to $\q$ (see remark \ref{remcanmodD}). 
\begin{remark}
This is the second notion of "reflex field". But the reader shouldn't get confused.
\end{remark}
Denote by $\Sigma$ the set of arithmetic subgroups $\Gamma$ of $GSp^{ad}(\q)^+$ which contain the image of a congruence subgroup of $GSp^{der}(\q)$.
The connected component of the identity $Sh^o$ of $Sh = Sh(GSp,\mathbb H^±_g)$ is then given by the inverse limit
$Sh^o = \varprojlim_{\Gamma \in \Sigma} \Gamma \backslash \mathbb H _g$ (cf. \ref{cshvD}). \\
Denote by $M^o=M^o(GSp,\mathbb H^±_g)$ the canonical model of $Sh^o$ in the sense of 2.7.10 \cite{Del2}, i.e. $M^o$ is defined over $\qbar$. For every $\Gamma \in \Sigma$ the space $\Gamma \backslash \mathbb H_g$ is an algebraic variety over $\komplex$ and $\Gamma \backslash M^o$ a model over $\overline \q$. \\
The field of rational functions $k(\Gamma \backslash M^o)$ on $\Gamma \backslash M^o$ is contained in the field of rational functions $k(\Gamma \backslash \mathbb H_g)$. Elements in the latter field correspond to meromorphic functions on $\sieg$ (now viewed as a complex analytic space) that are invariant under $\Gamma \in \Sigma$. \\
Following \cite{MiS} we call the field $ \mathcal {F} = \bigcup_{\Gamma \in \Sigma} k(\Gamma \backslash M^o)$ the field of \textbf{arithmetic automorphic functions} on $\sieg$. 

\subsubsection{About $Aut_\q(\mathcal F)$}

The (topological) group $\mathcal E$ defined by the extension (see 2.5.9 \cite{Del2}) \diagram{1 \ar[r] & \overline{GSp^{ad}(\q)^+} \ar[r] & \mathcal E \ar[r]^{\sigma \ \ \ \ \ \ \ } & Gal(\overline \q / \q) \ar[r] & 1} 
is acting continuously on $M^o$ (2.7.10 \cite{Del2}) and induces an action on $\mathcal F$ by \eq{ \label{action1} {}^{\alpha}\hspace{-0.6mm}f = \sigma(\alpha) \cdot (f \circ \alpha^{-1}) = (\sigma(\alpha)f) \circ (\sigma(\alpha)\alpha^{-1})} (see \cite{MiS} 3.2). This is meaningful because $f$ and $\alpha^{-1}$ are both defined over $\qbar$.
Using this action one can proove
\begin{theorem}[3.3 \cite{MiS}]
The map $\mathcal E \to Aut_\q(\mathcal F )$ given by (\ref{action1}) identifies $\mathcal E$ with an open subgroup of $Aut_\q(\mathcal F )$.
\end{theorem}

\subsection{Going down to $\q^{ab}$}
\label{goingdown}
We said that $M^o$ is defined over $\qbar$ but it is already defined over a subfield $k$ of $\q^{ab}$. More precisely $k$ is the fixed field of the kernel of the map $Gal(\qbar/\q)^{ab}\to \overline \pi _0 \pi(GSp)$ defined in 2.6.2.1 \cite{Del2}. 
Therefore the action of $\mathcal E$ on $M^o$ factors through the quotient $\overline{\mathcal E}$ of $\mathcal E$ defined by the following the commutative diagram with exact rows (see 4.2 and 4.12 \cite{MiS} or 2.5.3 \cite{Del2}) 
\diagramlab{diagtau}{
1 \ar[r]& \overline{GSp^{ad}(\q)^+} \ar[r] \ar[d]^{id} & \mathcal E \ar[r]^{\sigma \ \ \ \ \ } \ar[d]^{pr} & Gal(\qbar / \q) \ar[r] \ar[d]^{res} & 1 \\  
1 \ar[r] & \overline{GSp^{ad}(\q)^+} \ar[r] \ar[d]^{id} & \overline{\mathcal E} \ar[r]^{\sigma \ \ \ \ \ } \ar[d]^{ \tau} & Gal(k / \q) \ar[r] \ar[d]^{l} & 1 \\ 
1 \ar[r]& \overline{GSp^{ad}(\q)^+} \ar[r] & \frac{GSp(\adfin)}{\overline{C(\q)}} \ar[r] & \overline \pi _0 \pi(GSp) \ar[r] & 1,}
here $C$ denotes the center of $GSp$.
\begin{remark}
\label{remarkcenter}
1) We know that $C(\q)=\q^\times$ is discrete in $GSp(\adfin)$ (cf. \cite{hida}). \\
2) Because $l$ is injective, $\tau$ is injective as well and we can identify $\overline{\mathcal E}$ with an (open) subgroup of $\frac{GSp(\adfin)}{\q^\times}$. \\
3) $\overline {\mathcal E}$ is of course depending on $K$ but we suppress this dependence in our notation. 
\end{remark}

\subsubsection{The field $\mathcal M$ of arithmetic Modular functions}
\label{defmathcalM}
Let $f \in \mathcal F$ be a rational function, i.e. $f$ is a rational function on $\Gamma \backslash M^o$, for some $\Gamma \in \Sigma$. We call $f$ an \textbf{arithmetic Modular function} if it is rational over $k$, and meromorphic at the cusps (when viewed on the corresponding complex analytic space).  Compare this to 3.4 \cite{wei} or pp. 35 \cite{pedestrian}. 
\begin{definition}
The subfield of $\mathcal F$ generated by all arithmetic Modular functions is denoted by $\mathcal M$. Further we denote by $\mathcal M ^{cm}$ the subring of $\mathcal M$ of all arithmetic Modular functions which are defined in $x_{cm}$.
\end{definition}

The importance of $\mathcal M^{cm}$ for our purposes is explained (see 14.4 \cite{pedestrian} and 3.11 \cite{wei}) by
\begin{theorem}
\label{bigtheorem}
Let $A_{cm}$ denote the abelian variety corresponding to $x_{cm}$ (cf. \ref{intro2}). \\
Denote by \textbf{$K^{A_{cm}}$} the field extension of $K$ obtained by adjoining the field of definition of $A_{cm}$ and all of its torsion points. \\
Further denote by \textbf{$K^\mathcal M$} the field extension of $K$ obtained by adjoining the values $f(x_{cm})$, for $f \in \mathcal M ^{cm}$. \\
Finally denote by \textbf{$K^c$} the composition of $K$ with the fixed field of the image of the Verlagerungsmap $Ver : Gal(F^{ab}/F) \to Gal(E^{ab}/E)$, where $F$ is the maximal totally real subfield of $E$. Then we have the equality \eq{K^c=K^{A_{cm}}=K^\mathcal M.}
\end{theorem}

\begin{remark}
Notice we are not simply using the field of arithmetic automorphic functions as considered by Shimura, see 4.8 \cite{MiS} for his definition, because the exact size of the abelian extension obtained by using these functions is not clear (at least to the author). It is clear that the field of Shimura is contained in $K^\mathcal M$ and we guess that it should generate the same extension $K^c$ of $K$.
\end{remark}

\subsubsection{About $Aut_\q(\mathcal M)$}
\label{autmathcalM}
It is clear that $\mathcal M$ is closed under the action of $\overline {\mathcal E}$ (see 3.2 or 4.4 \cite{MiS}) and therefore we obtain a continuous map \eq{\overline {\mathcal E} \to Aut_\q(\mathcal M )} given like above by \eq{ \label{action2} {}^{\alpha}\hspace{-0.6mm}f = \sigma(\alpha) \cdot (f \circ \alpha^{-1}) = (\sigma(\alpha)f) \circ (\sigma(\alpha)\alpha^{-1}).}
In particular $GSp^{ad}(\q)^+$ is acting on $\mathcal M$ by \eq{ \label{actionQ} {}^{\alpha}\hspace{-0.6mm}f = f \circ \alpha ^{-1}.}

\subsection{The reciprocity law at $x_{cm}$}
\label{recilawgeneral}
Write \diagram{\mu_{cm} : \mathbb G _{m,\complex} \ar[r]^{} & S^E_\complex \ar[rr]^{ (h'_{cm})_\complex \  \ \ \ } & &\prod_{i=1}^r \mathcal T ^{E_{i,\complex}} \ar[r]^{\mathfrak i} & GSp_\complex} for the associated cocharacter of $h_{cm}$ (cf. \ref{canmodD}). From Lemma \ref{proph} 1) we know that $h_{cm} = \mathfrak i \circ \prod h_{\phi_i}$ and therefore $\mu_{cm} = \mathfrak i \circ \prod \mu_{\phi_i}$. Because $\mu_{\phi_i}$ is defined over $E_i^*$ the cocharacter $\mu_{cm}'=\prod \mu_{\phi_i}$ is defined over $\widetilde E = E_1^* \cdots E_r^* \subset E \subset K$. To simplify the notation set $\mathcal T = \prod^r_{i=1} \mathcal T^{E_i}$. 
Define the morphism \eq{\eta : T^K \to GSp} as composition of
\diagram{T^K \ar[rr]^{Res_{K / \q}(\mu'_{cm}) \ \ \  \ \ \ } & & Res_{K/ \q}(\mathcal T _K) \ar[rr]^{\ \ \ \ \ \ Nm_{K/\q}} && \mathcal T \ar[r]^{\mathfrak i \ \ } &GSp.} 
If we identify $\overline { \mathcal E}$ with an (open) subset of $\frac{GSp(\adfin)}{\q^\times}$ using $\tau$ (see (\ref{diagtau})), one can show (see 4.5 \cite{MiS} or 2.6.3 \cite{Del2}) that $\eta(\mathbb A):\mathbb A _{K}^\times \to GSp(\mathbb A)$ induces, by $\nu \mapsto \eta(\adfin)(\nu)$ mod $\q^\times$, a group homomorphism \eq{\label{overlineeta} \overline \eta : \mathbb A _{K}^\times \to \overline {\mathcal E}.}
If we denote by $[\nu] \in Gal(K^{ab} / K)$ the image of $\nu \in \mathbb A _{K}^\times$ under Artin's reciprocity map 
one can show (cf. 4.5 \cite{MiS}) that (cf. (\ref{diagtau})) \eq{ \label{cheapreci} \sigma(\overline \eta (\nu)) = [\nu]^{-1} | _k.}
\begin{remark}
The careful reader will ask why it is allowed to define $\eta$ using the extension $K$ of the field of definition of the cocharacter $\mu_{cm}$ given by $\widetilde E$, because our reference \cite{MiS} uses $\widetilde E$ to define $\eta$. The explanation for this is given by lemma \ref{lemmacomp} and standard class field theory.
\end{remark} 
Now we are able to state the \textbf{reciprocity law}
\begin{theorem}[see 4.6 and 4.10 \cite{MiS}]
\label{reclaw}
Let $\nu \in \mathbb A_K^\times$ and $f \in \mathcal M ^{cm}$. Then $f(x_{cm})$ is rational over $K^{ab}$. Further ${}^{\overline \eta (\nu)} \hspace{-1mm} f$ is defined in $x_{cm}$  and  \eq{\label{reci} {}^{\overline \eta (\nu)} \hspace{-1mm} f (x_{cm})=[\nu]^{-1}(f(x_{cm})).} 
\end{theorem}
\begin{proof}
We simply reproduce the argument given in the proof of Thm. 4.6 \cite{MiS}. The first assertion is clear by the definition of the canonical model (cf. \ref{canmodD}) and the other two assertions follow from the following calculation. \\
Regard the special point $x_{cm}$ as a point on the canonical model $[x_{cm},1] \in M^o$. The action of $\overline \eta (\nu)^{-1}$ is given by $\overline \eta (\nu)^{-1} [x_{cm},1] = \sigma(\overline \eta (\nu)^{-1})[x_{cm},\eta(\nu)]$ and further we know $[x_{cm},\eta(\nu)] = [\nu]^{-1}[x_{cm},1]$ (by (\ref{canmodspD})). Therefore we obtain 

\begin{align*}
{}^{\overline \eta (\nu)} \hspace{-1mm} f (x_{cm}) &= 
\sigma(\overline \eta (\nu)) \cdot  (f \circ \overline \eta (\nu)^{-1})([x_{cm},1]) \\ 
&= (\sigma(\overline \eta (\nu)f) \circ (\sigma(\overline \eta (\nu)) \overline \eta (\nu)^{-1}  ) ([x_{cm},1]) \\  
&= (\sigma(\overline \eta (\nu)f) \circ (\sigma(\overline \eta (\nu)) \sigma(\overline \eta (\nu))^{-1}) ([x_{cm},\eta(\nu)]) \\  
&\overset{\text{(\ref{cheapreci})}} = ([\nu]^{-1}|_k f)([\nu]^{-1}[x_{cm},1]) = [\nu]^{-1}(f([x_{cm},1])) \\ 
&= [\nu]^{-1}(f(x_{cm})).
\end{align*}
\end{proof}

The next observation is one of the key ingredients in our construction of the arithmetic subalgebra,
which we call therefore

\begin{proposition}
\label{mainprop}
The two maps of algebraic groups $\varphi$ and $\eta$ are equal.
\end{proposition}
\begin{proof}
This is an immediate corollary of Prop. \ref{propreflexnorm}, Lemma \ref{lemmaexplicitphi} , the compatibility properties of the norm map and the next simple Lemma.
\end{proof}

\begin{lemma}
\label{lemmacomp}
Let $(L,\phi)$ a CM type and $L'$ a finite extension of the reflex field $L^*$. Then the following diagram \diagram{
T^{L'} \ar[rr]^{Res_{L' / \q}{\mu_{L'}} \ \ \ \ \ \ } \ar[d]^{N_{L' / L^*}} && Res_{L' / \q}(T^L_{L'}) \ar[rr]^{ \ \ \ \ \ Nm_{L' / \q}} & & T^L \\
T^{L^*} \ar[rr]^{Res_{L^* / \q}{\mu_{L^*}} \ \ \ \ \ \ } & & Res_{L^* / \q}(T^L_{L^*}) \ar[rru]_{ \ \ \ \ \ \ Nm_{L^* / \q}} & &} 
is commutative.
\end{lemma}

Now having all the number-theoretic ingredients we need in hand, we can move on to the "operator-theoretic" part of this paper.

%%%%%%%%%%%%%%%%%%%%%%%%%%%%%%%%%%%%%%%%%%%%%
%%%%%%%%        REVIEW OF HA PAUGAM
%%%%%%%%%%%%%%%%%%%%%%%%%%%%%%%%%%%%%%%%%%%%%

\section{On Bost-Connes-Marcolli systems}
\label{hapabcm}
We review very briefly the general construction of $C^*$-dynamical systems, named Bost-Connes-Marcolli systems, as given in \cite{HaPa}. 

\subsection{BCM pairs}
\label{bcmpairs}
A BCM pair $(\mathcal D, \mathcal  L)$ is a pair consisting of a BCM datum $\mathcal D = (G,X,V,M)$ together with a level structure $\mathcal L = (L,\Gamma,\Gamma_M)$ of $\mathcal D$. \\
A BCM datum is a Shimura datum $(G,X)$ together with an enveloping algebraic semigroup $M$ and a faithful representation $\phi:G \to GL(V)$ such that $\phi(G) \subset M \subset End(V)$. Here $V$ denotes a $\q$-vector space of finite dimension. \\
A level structure $\mathcal L$ of $\mathcal D$ consists of a lattice $L \subset V$, a compact open subgroup $\Gamma \subset G(\adfin)$ and a compact open semigroup $\Gamma_M \subset M(\adfin)$ such that $\phi(\Gamma) \subset \Gamma_M$ and $\Gamma_M$ stabilizes $L\otimes_\ganz \hat \ganz$. 
\begin{remark}
In 3.1 \cite{HaPa} a more general notion of Shimura datum is allowed than ours given in the Appendix.
\end{remark}
To every BCM datum $\mathcal D$ and lattice $L \subset V$ one can associate the following so called \textbf{maximal level structure} to obtain a BCM pair by setting $\Gamma_{M} = M(\adfin) \cap End(L\otimes_\ganz \widehat {\mathbb Z})$ and $\Gamma_{} = \phi^{-1}(\Gamma_M^\times)$. \\ \\
The level structure $\mathcal L$ is called \textbf{fine} if $\Gamma$ is acting freely on $G(\q) \backslash (X \times G(\adfin))$.   
\begin{remark}
For the definition of the topology of $G(\adfin)$ and $M(\adfin)$ we refer the reader to \cite{plat}. Especially one can show that $\phi(\adfin) : G(\adfin) \to M(\adfin)$ is a continuous map (cf. lemma 5.2 \cite{plat}).
\end{remark}

\subsection{Quotient maps attached to BCM pairs}

Let $(\mathcal D,\mathcal L)$ a BCM pair.

\subsubsection{The BCM groupoid}

There is a partially defined action of $G(\adfin)$ on the direct product $\Gamma_M \times Sh(G,X)$ given by \eq{g (\rho,[z,l]) = (g\rho,[z,lg^{-1}]),}
where we suppressed the mophism $\phi$. Using this the BCM groupoid $U$ is the topological groupoid (using the notation given in (\ref{notlln})) defined by \eq{U = G(\adfin) \boxplus (\Gamma_M \times Sh(G,X)).}
There is an action of the group $\Gamma \times \Gamma$ on $U$ given by \eq{\label{actgeneralgroupoid}(\gamma_1,\gamma_2)(g,\rho,[z,l]) (\gamma_1 g \gamma_2^{-1},\gamma_2 \rho,[z,l\gamma_2^{-1}]).}
We denote the quotient by $Z=(\Gamma \times \Gamma) \backslash U$ and obtain a natural quotient map \diagram{U \ar[r] & Z.}
We denote elements in $Z$ by $[g,\rho,[z,l]]$.
\begin{remark}
In general the quotient $Z$ is not a groupoid anymore, see 4.2.1 \cite{HaPa}.
\end{remark}

\subsubsection{The positive BCM groupoid}
Assume that the Shimura datum $(G,X)$ of our BCM pair satisfies $(SV5)$ (cf. \ref{datumD}). Moreover we choose a connected component $X^+$ of $X$ and set $G(\q)^+ = G(\q) \cap G(\real)^+$, where $G(\real)^+$ denotes the connected component of the identity of $G(\real)$. Then $G(\q)^+$ is acting naturally on $X^+$, because $X^+$ can be regarded as a $G(\real)^+$-conjugacy class (see \ref{cshvD}). Now we can consider the positive (BCM) groupoid $U^+$ which is the topological groupoid given by \eq{U^+ = G(\q)^+ \boxplus (\Gamma_M \times X^+).}
If we set $\Gamma^+ = \Gamma \cap G(\q)^+$ we see further that $\Gamma^+ \times \Gamma^+$ is acting on $U^+$ by \eq{\label{actplus}(\gamma_1,\gamma_2)(g,\rho,z) = (\gamma_1 g \gamma_2^{-1},\gamma_2 \rho,\gamma_2z).}
We denote the quotient by $Z^+ = (\Gamma^+ \times \Gamma^+) \backslash U^+$
and obtain another quotient map \diagram{U^+ \ar[r]& Z^+.} 
There is a natural equivariant morphism of topological groupoids \diagramlab{naturalmorph}{U^+ \ar[r] &  U} given by $(g,\rho,z) \mapsto (g,\rho,[z,1])$, inducing a commutative diagram \diagramlab{commdiagplus}{U \ar[r] & Z \\ U^+ \ar[u] \ar[r] & Z^+  \ar[u]}

The following criterion given in 5.1 \cite{HaPa} will be crucial for our approach
\begin{criterion}
\label{crit}
If the natural map $G(\q) \cap \Gamma \to G(\q)/G(\q)^+$ is surjective and $| G(\q)\backslash G(\adfin) / \Gamma  |  = 1$, then the natural morphism (\ref{naturalmorph}) induces a homeomorphism of topological spaces \diagram{Z^+ \ar[r] & Z.}
\end{criterion}

\begin{remark}
\label{decompositionandexplicitinverse}
1) In other words the two conditions of the criterion simply mean that we have a decomposition of the form $G(\adfin)= G(\q)^+ \cdot \Gamma$. \\
2) The inverse $Z \longrightarrow Z^+$ of the above homeomorphism is given explicitly as follows. By using the first remark we can write every $l \in G(\adfin)$ as a product $l = \alpha \beta$ with $\alpha \in G(\q)^+$ and $\beta \in \Gamma$ (this decomposition is unique up to an element in $\Gamma^+ = \Gamma \cap G(\q)^+$). In particular every element $[g,\rho,[z,l]] \in Z$ can be written as $[g,\rho,[z,l]] = [g\beta^{-1},\beta\rho,[\alpha^{-1}z,1]]$ and under the inverse of the above homeomorphism this element is sent to $[g\beta^{-1},\beta\rho,\alpha^{-1}z] \in Z^+$.
\end{remark}

\subsubsection{The adjoint BCM algebra}

Let us denote by $C$ the center of $G$ and \textbf{assume} further that $\phi(\overline{C(\q)})$ is a normal subsemigroup of $M(\adfin)$. The adjoint group $G^{ad}$ of $G$ is the quotient of $G$ by its center $C$ (in the sense of algebraic groups, see \cite{water}). Let us define the semigroup $\Gamma^{ad}_M$ to be the quotient of $\Gamma_M$ by the normal subsemigroup $\phi(\overline{C(\q)}) \cap \Gamma_M$ and remember that $X^+$ can be naturally regarded as a $G^{ad}(\real)^+$-conjugacy class (see \ref{cshvD}). With this in hand we define the adjoint (BCM) groupoid $U^{ad}$ to be the topological groupoid \eq{U^{ad} = G^{ad}(\q)^+ \boxplus (\Gamma^{ad}_M \times X^+).}
It is known that the projection $G \longrightarrow G^{ad}$ induces a surjective group homomorphism $\pi^{ad} : G(\q)^+ \longrightarrow G^{ad}(\q)^+$ (see 5.1 \cite{isv}). Setting $\Gamma^{ad} = \pi^{ad}(\Gamma^+)$ we see immediately that $\Gamma^{ad}\times \Gamma^{ad}$ is acting on $U^{ad}$ exactly as in \ref{actplus}. We obtain yet another quotient map \diagram{U^{ad} \ar[r] &  Z^{ad}.}
Using the two projections $\Gamma_M \to \Gamma^{ad}_M$ and $\pi^{ad}$ there is by construction an obvious equivariant morphism of topological groupoids \eq{\label{naturalmorph2}U^+ \longrightarrow U^{ad}} which induces (together with (\ref{commdiagplus})) a commutative diagram \diagram{U \ar[r] & Z \\ U^+ \ar[u] \ar[r] \ar[d] & Z^+  \ar[u] \ar[d] \\ U^{ad} \ar[r] & Z^{ad}. }

\subsection{BCM algebras and systems}
\label{bcmalgebras}
Let $(\mathcal D,\mathcal L)$ be a BCM pair. The \textbf{BCM algebra} $H = H_{\mathcal D,\mathcal L}$ is defined to be the set of compactly supported, continuous function on the quotient $Z = (\Gamma \times \Gamma) \backslash U$ of the BCM groupoid $U$, i.e. \eq{H=C_c(Z).}
By viewing functions in $H$ as $\Gamma \times \Gamma$-invariant functions on the groupoid $U$, we can equip $H$ with the structure of a $*$-algebra by using the usual convolution and involution on $U$ (like in the construction of groupoid $C^*$-algbras). We refer to 4.3.2 \cite{HaPa} for the details. After completing $H$ in a suitable norm we obtain a $C^*$-algebra $A$ (see 6.2 \cite{HaPa}). Further there is a time evolution $(\sigma_t)_{t \in \real}$ on $H$ resp. $A$ so that we end up with the \textbf{BCM system} $\mathcal A = \mathcal A _{\mathcal D,\mathcal L}$ given by the $C^*$-dynamical system \eq{\mathcal A = (A,(\sigma_t)_{t\in \real})}
associated with the BCM pair $(\mathcal L,\mathcal D)$. For the general definition of the time evolution we refer to 4.4 \cite{HaPa}. We will state the \textbf{time evolution} $(\sigma_t)_{t \in \real}$ only in the case of our BC-type systems $\mathcal A _K$ (cf. \ref{timeevolK}).

\begin{remark}
In complete analogy one might construct a positive respectively adjoint BCM system, but we don't need this.
\end{remark}

\subsection{On Symmetries of BCM algebras}

In section 4.5 \cite{HaPa} the authors define symmetries of BCM algebras but for our purpose we need to deviate from their definition in order to be in accordance with the definition of symmetries for BC-type systems given in \cite{Sergey}.  \\ \\
Let $(\mathcal D, \mathcal L)$ be a BCM pair with fine level structure (see \ref{bcmpairs}) and recall (\ref{shvD}) that there is a natural right action of $G(\adfin)$ on the Shimura variety $Sh(G,X)$ which is denoted by $m [z,l] = [z,lm]$. Define the subgroup $G_\Gamma(\adfin) = \{ g \in G(\adfin) \ | \ g \gamma = \gamma g \ \forall \gamma \in \Gamma \}$. Further, if we denote by $C$ the center of $G$, the group $C(\real)$ is acting on $Sh(G,X)$ by $c[z,l] = [cz,l]$. We end up with a right action of $G_\Gamma(\adfin) \times C(\real)$ as symmetries on the BCM algebra $H_{\mathcal D, \mathcal L}$ given on a function $f \in C_c(Z)$ by \eq{^{(m,c)} \hspace{-1mm}f(g,\rho,[z,l]) = f(g,\rho,[cz,lm]).}
\begin{remark}
If $G(\adfin)$ is a commutative group and we have a decomposition $G(\adfin) = G(\q) \cdot \Gamma$ then it is immediate that our symmetries agree with the ones defined in 4.5 \cite{HaPa}.
\end{remark}

\section{Two BCM pairs and a map}

In this section we will apply the constructions from the last section to our Shimura data $\mathcal S_K = (T^K,X_K,h_K)$ and $\mathcal S_\siegel = (GSp(V_E,\psi_E),\mathbb H ^±_g,h_{cm})$ from section \ref{shimdata}
and show how the two resulting systems can be related.

\subsection{Costume I: $(\mathcal D_K,\mathcal L_K)$ and $\mathcal A_K$}
\label{mathcalAK}
This section is valid for an arbitrary number field $K$.
\subsubsection{The BCM groupoid $(\mathcal D_K,\mathcal L_K)$}
Let us recall the BCM pair $(\mathcal D _K,\mathcal L _K)$ from 5.5 \cite{HaPa} attached to $\mathcal S_K$. It is given by \eq{(\mathcal D _K,\mathcal L _K) = ((\mathcal S_K,K,M^K),(\mathcal O_K,\widehat{\mathcal O}_K^\times, \widehat{\mathcal O}_K) ),}
where the algebraic semigroup $M^K$ is represented by the functor which assigns to a $\q$-algebra $R$ the semigroup of $\q$-algebra homomorphisms $Hom(K[X],K \otimes_\q R)$. By definition we have that $M^K(R)^\times = T^K(R)$ for every $\q$-algebra $R$, which gives an embedding $\phi : T^K \to M^K$. (It will be convenient to set $\Gamma_K = \widehat{\mathcal O}_K^\times$.)

\subsubsection{The quotient map $U_K \to Z_K$}

The corresponding BCM groupoid, denoted by $U_K$, is given by \eq{U_K = T^K(\adfin)\boxplus(\widehat{\mathcal O}_K \times Sh(\mathcal S _K))}
with $\Gamma_K^2 = \widehat{\mathcal O}^\times_K \times \widehat{\mathcal O}^\times_K$ acting as in (\ref{actgeneralgroupoid}). We denote the quotient of this action by \eq{Z_K = \Gamma_K^2 \backslash U_K.}

\subsubsection{The time evolution}
\label{timeevolK}
Following 7.1 \cite{HaPa} the time evolution $(\sigma_t)_{t\in \real}$ on the BCM algebra $H_K=C_c(Z_K)$ is given as follows. Denote by $\mathfrak N = \mathfrak N_{K/\q}:\mathbb A _{K,f}^\times \to \real$ the usual idele norm. Let $f \in H_K$ be a function, then we have \al{\sigma_t(f)(g,\rho,[z,l]) = \mathfrak N (g)^{it} f(g,\rho,[z,l]).}

\subsubsection{On symmetries}
\label{symmbcmK}

First, from 2.2.3 \cite{Del2} we know that there is an isomorphism between $Sh = Sh(T^K,X_K,h_K)$ and $\pi_0(C_K)$. By class field theory the latter space $\pi_0(C_K) = C_K / D_K$ is identified with the Galois group $Gal(K^{ab}/ K)$ of the maximal abelian extension of $K$ using the Artin reciprocity homomorphism. Under this identification the natural action of $T^K(\adfin) = \mathbb A _{K,f}^\times$ on $Sh$ corresponds simply to the Artin reciprocity map, i.e. if $\nu$ is a finite idele in $T(\adfin)$ and $\omega_1 = [g,l] \in Sh$ corresponds to the identity in $Gal(K^{ab}/K)$, then $\nu \omega_1 = [g,l \nu]$ corresponds to the image $[\nu]$ in $Gal(K^{ab} / K)$ of $\nu$ under Artin's reciprocity map. \\ \\
Now, because $T^K$ is commutative we see that $C(\adfin)\times C(\real) = T^K(\mathbb A) = \mathbb A _K^\times$ is acting by symmetries on $H_K$. By what we just said this action is simply given by the natural map $\mathbb A^\times _K \to \pi_0(C_K) = C_K / D_K$ so that we obtain (tautologically) the desired action of $C_K / D_K \cong Gal(K^{ab} / K)$ on $H_K$.

\begin{remark}
The reader should notice that in the case of an imaginary quadratic number field $K$ our symmetries do not agree with the symmetries defined \cite{cmr1} (except when the class number of $K$ is equal to one, where the two definitions agree). For a short discussion on this matter we refer the reader to remark \ref{finalremark}.
\end{remark}

\subsubsection{About extremal $KMS_\infty$-states of $\mathcal A_K$}
\label{extrkms}
We refer to pp. 445 \cite{conmar} or \cite{brat2} for the notion of extremal $KMS_\infty$-states. \\ \\
Let $\mathcal A _K = (A_K,(\sigma_t)_{t\in \real})$ denote the corresponding BCM system (cf. \ref{bcmalgebras}). In Theorem 2.1 (vi) \cite{Sergey} it is shown that the set $\mathcal E _\infty$ of extremal $KMS_\infty$-states of $\mathcal A _K$ is indexed by the set $Sh = Sh(T^K,X_K,h_K)$ and the extremal $KMS_\infty$-state $\varrho_\omega$ associated with $\omega \in Sh$ is given on a function $f \in H_K$ by evaluation, namely \eq{\varrho_\omega(f) = f(1,1,\omega).}

\begin{remark}
1) It follows immediately that the symmetry group $C_K / D_K$ is acting free and transitively on the set of extremal $KMS_\infty$-states. \\
2) Using the definition of symmetries given in \cite{HaPa} we obtain another action of $C_K / D_K$ on $H_K$ but in general this will not induce a free and transitive action on the extremal $KMS_\infty$-states. In fact the two different actions of $C_K / D_K$ on $H_K$ given in \cite{Sergey} resp. \cite{HaPa} are equivalent if and only if the class number $h_K$ of $K$ is equal to $1$. This follows directly from the fact that if (and only if) $h_K = 1$ then already $\widehat {\mathcal O} _K^\times$ surjects onto $Gal(K^{ab}/K)$. 
\end{remark}

All put together we have

\begin{theorem}[\cite{HaPa} and \cite{Sergey}]
Let $K$ be an arbitrary number field.
Then the BCM system $\mathcal A_K = (A_K,(\sigma_t)_{t\in \real})$ is a BC-type system (cf. page 2).
\end{theorem}

\subsection{Costume II: $(\mathcal D _\siegel,\mathcal L _\siegel)$}

\subsubsection{The BCM groupoid $(\mathcal D _\siegel,\mathcal L _\siegel)$}
\label{costumeII}
Recall the construction of the symplectic vector space $(V_E,\psi_E)$ (see \ref{construction}). We still have some freedom in specifying the totally imaginary generators $\xi_i$ of the (primitive) CM fields $E_i$, which in turn define the symplectic form $\psi_E$ (cf. (\ref{defpsi})). Let us denote by $L_E$ the lattice $L_E = \oplus_{i=1}^r \mathcal O _{E_i}\subset V_E= \oplus_{i=1}^r E_i$. We now fix generators $\xi_i$ according to 

\begin{lemma}
\label{lemmaintegralbasis}
For each $i \in \{1,..,r\}$ there exists a totally imaginary generator $\xi_i \in E_i$, such that the associated symplectic vector space $(V_E,\psi_E)$ is integral with respect to $L_E$, i.e. there exists a symplectic basis $\{e_j\}$ for $(V_E,\psi_E)$ such that $e_j \in L_E$, for each $j$.
\end{lemma} 
\begin{proof}
For each $i$, choose any totally imaginary generator $\widetilde{\xi}_i \in E_i$ and regard the associated symplectic form $\widetilde {\psi}_E$ on $V_E$ (see (\ref{defpsi})).
It is known that there exists a symplectic basis $\{\widetilde e _j\}$ for $(V_E, \widetilde \psi_E)$. Now for each $j$ there exists a $q_j \in \mathbb N$ such that $e_j = q_j \widetilde e_j \in L_E $, because $E_i = \mathcal O _{E_i}\otimes_\ganz \q$. Set $q = \prod_j q_j$ and define the symplectic form $\psi_E$ on $V_E$ by using the totally imaginary generators $\xi _i= q^{-2}\widetilde \xi _i \in E_i$, for every $i$. By construction it is now clear that $\{e_j\}$ is an integral symplectic basis for $(V_E,\psi_E)$.
\end{proof}

Now having fixed our Shimura datum $\mathcal S _\siegel =(GSp(V_E,\psi_E),\mathbb H_g^±,h_{cm})$ we define the BCM pair $(\mathcal D_\siegel,\mathcal L_\siegel)$ equipped with the maximal level structure (cf. \ref{bcmpairs}) with respect to the lattice $L_E$ by \eq{(\mathcal D_\siegel,\mathcal L_\siegel) = ( (\mathcal S _\siegel,V_E,MSp),(L_E,\Gamma_{\siegel}, \Gamma_{\siegel,M}) ),}
where the algebraic semigroup $MSp = MSp(V_E,\psi_E)$ is represented by the functor which assigns to a $\q$-algebra $R$ the semigroup 
\alst{&MSp(R) = \\ &\{ f\in End_R(V_E \otimes_\q R ) \ | \ \exists \ \nu(f) \in R \text{ : } \psi_{E,R}(f(x),f(y)) = \nu(f) \psi_{E,R}(x,y) \ \forall x,y \}.}
It is clear by definition (compare \ref{exA}) that $MSp(R)^\times = GSp(R)$ which defines a natural injection $\phi : GSp \to MSp$.
\subsubsection{Some quotient maps}
\label{somequotientmaps}
We denote the corresponding BCM groupoid by \eq{U_\siegel = GSp(\adfin)\boxplus(\Gamma_{\siegel,M} \times Sh(\mathcal S _\siegel)),}
where the group $\Gamma_\siegel^2 = \Gamma_\siegel \times \Gamma_\siegel$ is acting as usual. We denote the quotient of $U_\siegel$ by this action by \eq{Z_\siegel = \Gamma_\siegel^2 \backslash U_\siegel.}
Thanks to \ref{exD} and Remark \ref{remarkcenter} 1) we are allowed to consider the positive and adjoint BCM groupoid, which we denote by $U_\siegel^+$ and $U_\siegel^{ad}$ respectively. The corresponding quotients are denoted analogously by $Z_\siegel^+ = (\Gamma_\siegel^+)^2 \backslash U_\siegel^+$ and $Z_\siegel^{ad}=(\Gamma_\siegel^{ad})^2 \backslash U_\siegel^{ad}$.

\subsection{The map $\Theta : Z_K \to Z_\siegel^{ad}$}
\label{thebcmmap}
The aim in this section is to construct a continuous map $\Theta : Z_K \longrightarrow Z_\siegel^{ad}$.

\subsubsection{Relating $Z_K$ and $Z_\siegel$}
\label{continuation}
Recall that the morphism of Shimura data $\varphi : \mathcal S_K \to \mathcal S_\siegel$ constructed in \ref{themapphi} is a morphism of algebraic groups $\varphi : T^K \to GSp$, inducing a morphism $Sh(\varphi) : Sh(\mathcal S_K)\to Sh(\mathcal S_\siegel)$ of Shimura varieties. Moreover there is a natural continuation of $\varphi$ to a morphism of algebraic semigroups $M(\varphi) : M^K \to MSp$ due to the following. We know that $\varphi:T^K\to GSp$ can be expressed in terms of reflex norms (see lemma \ref{lemmaexplicitphi} and \ref{reflexB}), which are given by determinants and this definition makes still sense if we replace $T^K$ and $GSp$ by their enveloping semigroups $M^K$ and $MSp$, respectively. Now we can define an equivariant morphism of topological groupoids \al{\Omega : U_K  \to U_\siegel} by \al{ \label{mapbcmgroupoid} (g,m,z) \in U_K \mapsto (\varphi(\adfin)(g),M(\varphi)(\adfin)(z),Sh(\varphi)(z)) \in U_{\siegel}.} 
To show the equivariance of $\Omega$ use $\varphi(\adfin)(\Gamma_K)\subset \Gamma_\siegel$ and the equivariance of $Sh(\varphi)$ (cf. \ref{shvD}). We obtain a continuous map \al{\overline{\Omega} : Z_K \to Z_\siegel.}

\subsubsection{Relating $Z_\siegel$ and $Z_\siegel^{ad}$} 

In order to relate $Z_\siegel$ and $Z_\siegel^{ad}$ we will show that we are allowed to apply Criterion \ref{crit} by proving the following two lemmata. 
\begin{lemma}
\label{lemmadecompplus}
We have $GSp(V_E,\psi_E)(\adfin) = GSp(V_E,\psi_E)(\Q) \cdot \Gamma_\siegel .$
\end{lemma}
\begin{proof}
Let $\{ e_{j}\}$ be an integral symplectic basis of $V_E$ with respect to $L_E$ (cf. Lemma \ref{lemmaintegralbasis}). Each $f \in GSp(\adfin)$ is $\adfin$-linear and therefore determined by the 
values on $e_j \otimes 1 \otimes 1 \in L_E \otimes_\ganz \Q \otimes_\ganz \widehat \ganz$ given by \al{f(e_j\otimes 1 \otimes 1) = \sum_k a_{k,j} \otimes b_{k,j} \otimes c_{k,j} \in L_E \otimes_\ganz \Q \otimes_\ganz \widehat \ganz.}
Let $d_{k,j} \in \mathbb N$ be the denominator of $b_{k,j}$ and define $c(f) = \prod_{k,j} d_{k,j} \in \mathbb N.$ \\
Now observe that the map $$M_f : e_i \mapsto c(f) e_i$$ is a map in $GSp(\Q) \subset GSp(\adfin)$
i.e. $M_f$ is compatible with the symplectic struture $\psi_E$.
Obviously $$M_f \circ f \in \Gamma_{\siegel}$$
and thus we obtain the desired decomposition $$f = M_f^{-1} \circ (M_f \circ f) \in GSp(\Q) \cdot \Gamma_\siegel.$$
\end{proof}

\begin{lemma}
The map $GSp(\q) \cap \Gamma_\siegel \to GSp(\Q)/GSp(\Q)^+$ is surjective.
\end{lemma}
\begin{proof}
We know that $GSp(\reell)^+ = \{f\in GSp(\reell) | \nu(f) > 0\}$ (see \ref{exA}). \\
From this we get $GSp(\Q)^+ = GSp(\q) \cap GSp(\real)^+ =\{ f\in GSp(\Q) | \nu(f) > 0 \}.$ \\
Let $f$ be an element in $GSp(\Q)$. If we define $M_f$ exactly like in the proof above, we see that $M_f \in GSp(\q)^+$ and can conclude $M_f \circ f \in GSp(\q) \cap \Gamma_\siegel.$
\end{proof}

Thus we see that the natural morphism $Z^+_\siegel \to Z_\siegel$ (see \ref{naturalmorph}) is a homeomorphism so that we can invert this map and compose it with the natural map $Z^+_\siegel \to Z^{ad}_\siegel$ (cf. \ref{naturalmorph2}) to obtain a continuous map \diagram{Z_\siegel \ar[r] & Z^{ad}_\siegel.}
Finally if we compose the last map with $\overline \Omega$ from above we obtain a continuous morphism denoted \diagram{\Theta : Z_K \ar[r] & Z_\siegel^{ad}.}
One crucial property of $\Theta$ is that every element $z \in Z_K$ will be sent to an element of the form \eq{\label{imagetheta} \Theta(z)=[g\beta^{-1},\beta \rho, \alpha^{-1}x_{cm}] \in Z^{ad}_\siegel,} where $g \in G^{ad}(\q)^+$, $\rho \in \Gamma^{ad}_{\siegel,M}$, $\alpha \in G^{ad}(\q)$ and $\beta \in \Gamma^{ad}_\siegel$, such that $\alpha \beta \in \pi^{ad}(\varphi(T^K(\adfin))) \subset G^{ad}(\adfin)$. \\ \\
Now we are ready for the

\section{Construction of our arithmetic subalgebra}
\label{constrarith}
The idea of the construction to follow goes back to \cite{cmr1} and \cite{cmr2}. \\ \\
We constructed the ring $\mathcal M ^{cm}$ of arithmetic Modular functions on $\mathbb H_g$ that are defined in $x_{cm}$ (see \ref{defmathcalM}). Further the group $\overline{\mathcal E}$ acts by automophisms on $\mathcal M ^{cm}$ according to \ref{autmathcalM}. (Recall that we use the notation ${}^{\alpha}\hspace{-0.5mm}f$ to denote the action of an automorphism $\alpha$ on a function $f \in \mathcal M^{cm}$.) Thanks to the embeddings (cf. \ref{goingdown}) \diagram{\overline{\mathcal E} \ar[r]& \frac{GSp(\adfin)}{\q^\times} \ar[r] & \frac{MSp(\adfin)}{\q^\times} & \Gamma_{\siegel,M}^{ad} \ar[l]} 
the intersection $\overline {\mathcal E } \cap \Gamma_{\siegel,M}^{ad} $ is meaningful and thus we can define, for each $f \in \mathcal M ^{cm}$, a function $\widetilde f$ on $U_\siegel^{ad}$ by \eq{ \label{defarithel} \widetilde f (g,\rho,z) =   \left\{ \begin{array}{cc} ^{\rho}\hspace{-0.5mm}f(z), & \text{if } \rho \in \overline {\mathcal E} \cap \Gamma_{\siegel,M}^{ad}  \\ 0 & \text{else.} \end{array} \right.}
By construction $\widetilde f$ is invariant under the action of $(\gamma_1,\gamma_2) \in \Gamma_\siegel^{ad}\times \Gamma_\siegel^{ad}$ because
\eq{\widetilde f (\gamma_1 g \gamma_2^{-1},\gamma_2 m, \gamma_2 z) \overset{\text{def}}{=} {}^{\gamma_2 m}\hspace{-0.6mm}f(\gamma_2 z)\overset{(\ref{actionQ})}{=} {}^{\gamma_2^{-1}\gamma_2 m}\hspace{-0.6mm}f(z) = \widetilde f(g,m,z).}
Therefore we can regard $\widetilde f$ as function on the quotient $Z_\siegel^{ad}=(\Gamma^{ad}_\siegel)^2 \backslash U_\siegel^{ad}$. \\ \\
We set $W_K = \widehat {\mathcal O}_K^\times \times \widehat {\mathcal O}_K \times Sh(T^K,X_K)$, which is a compact and clopen subset of $U_K$  and invariant under the action of $\Gamma_K^2$, i.e. $\Gamma_K^2 \cdot W_K \subset W_K$. With these preliminaries we have the following

\begin{proposition}
\label{pullbackprop}
Let $f$ be a function in $\mathcal M^{cm}$, then $\widetilde f \circ \Theta$ is contained in $C_c(Z_K)$, i.e. \eq{\widetilde f \circ \Theta \in H_K \subset A_K.}
\end{proposition}
\begin{proof}
As we already remarked in (\ref{imagetheta}) the image of an element $z \in Z_K$ under $\Theta$ is of the form $[g\beta^{-1},\beta \rho, \alpha^{-1}x_{cm}] \in Z^{ad}_\siegel$. Therefore $f_K = \widetilde f \circ \Theta$ is continous, because $\Theta$ is continuous, the action of $\overline{\mathcal E}$ is continuous and does not produce singularities at special points (see Theorem \ref{reclaw}). \\
Let us now regard $f_K$ as $\Gamma_K^2$-invariant function on $U_K$. Thanks to $\overline {\mathcal E} \subset \frac{GSp(\adfin)}{\q^\times}$ and (\ref{defarithel}) we see that the support of our function $f_K$ is contained in the clopen subset $\widehat {\mathcal O}_K^\times \times \widehat {\mathcal O}_K^\times \times Sh(T^K,X_K) \subset U_K$. Using the compact subset $W_K \subset U_K$, the next easy lemma finishes the proof.
\end{proof}

\begin{lemma}
Let $G$ be a topological group, $X$ be a topological $G$-space and $Y \subset X$ a compact, clopen subset such that $G Y \subset Y$. If we have a continuous, $G$-invariant function $f \in C^G(X)$ then $ f |_Y \in C_c(G \backslash X)$.
\end{lemma}

Now we can define our arithmetic subalgebra of $\mathcal A _K=(A_K,\sigma_t)$ (cf. \ref{mathcalAK}).

\begin{definition}
Denote by $A _K^{arith}$ the $K$-rational subalgebra of $A_K$ generated by the set of functions $\{ \widetilde f \circ \Theta \ | \ f \in \mathcal M^{cm}\}$.
\end{definition}

\subsection{Proof of Theorem \ref{theorem}}
\label{proof}
Let $f \in \mathcal M ^{cm}$ and denote $f_K = \widetilde f \circ \Theta \in A_K^{arith}$. Further denote by $\varrho_\omega$ the extremal $KMS_\infty$-state of $\mathcal A _K$ corresponding to $\omega \in Sh = Sh(T^K,X_K,h_K)$ (see \ref{extrkms}). Recall the isomorphism $Gal(K^{ab}/ K) \cong \pi_0(C_K) = Sh$ given by Artin reciprocity. Considered as element in $Gal(K^{ab}/K)$ we write $[\omega]$ for $\omega$.

\subsubsection*{Property vi)}
\label{proof2}
Let $\nu \in \mathbb A_K^\times$ be a symmetry of $\mathcal A _K$ (see \ref{symmbcmK}). Thanks to Lemma \ref{lemmaintegralbasis} and \ref{lemmadecompplus}, we can write $\varphi(\adfin)(\nu) = \alpha \beta \in GSp(\adfin)$ with $\alpha \in GSp(\q)^+$ and $\beta \in \Gamma_\siegel$. By $\overline \alpha$ resp. $\overline \beta$ we will denote their images in $GSp^{ad}(\q)^+$ resp. $\Gamma_\siegel^{ad}$ under the map $\pi^{ad} \circ \varphi(\adfin)$.  Moreover we denote the image of $\nu$ under Artin reciprocity by $[\nu] \in Gal(K^{ab}/ K)$.\\
The action of the symmetries on the extremal $KMS_\infty$-states is given by pull back and because this action is free and transitive it is enough to restrict to the case of the extremal $KMS_\infty$-state $\varrho_1$ corresponding to the identity in $Gal(K^{ab}/K)$. \\
Using Proposition \ref{mainprop} and the reciprocity law (\ref{reci}) we can calculate the action of $\nu$ on $\varrho_\omega(f_K)$ as follows

\begin{align*}
{}^{\nu}\hspace{-0.6mm}\varrho_1(f_K) &\overset{\text{def}}{=}
\varrho_1(^{\nu}\hspace{-0.6mm}f_K) 
\overset{(\ref{imagetheta})}=\widetilde f (\overline \beta ^{-1},\overline \beta,\overline \alpha^{-1}x_{cm}) 
\overset{\text{\ref{mainprop}}}{=} {}^{\overline \eta(\beta)}\hspace{-0.6mm}f(\overline\alpha^{-1}x_{cm}) \\
&\overset{(\ref{actionQ})}={}^{\overline \eta(\alpha)\overline \eta(\beta)}\hspace{-0.6mm}f(x_{cm}) 
\overset{}={}^{\overline \eta(\nu)}\hspace{-0.6mm}f(x_{cm}) 
\overset{\text{(\ref{reci})}}{=} [\nu]^{-1}(f(x_{cm})) \\
&=[\nu]^{-1}(\varrho_1(f_K))
\end{align*}

This is precisely the intertwining property we wanted to prove.

\subsubsection*{Property v)}
\label{proof1}
Using the notation from above, we conclude immediately that by construction and Theorem \ref{bigtheorem} we have 
\eq{\varrho_1(f_K) = \widetilde f (1,1,x_{xm}) = f(x_{cm}) \in K^{ab}}
and the above calculation shows further that 
 \eq{\varrho_\omega(f_K) = {}^{[\omega]} \hspace{-0.6mm} \varrho_1(f_K) = [\omega]^{-1}(f(x_{cm})) \in K^{ab}}
finishing our proof.

\begin{remark}
\label{finalremark}
We want to conclude our paper with a short discussion comparing our construction with the original construction of Connes, Marcolli and Ramachdran (see \cite{cmr1}) in the case of an imaginary quadratic field $K$. Apart from the fact that we are not dealing with the "$K$-lattice" picture as done in \cite{cmr1} the main difference lies in the different definitions of symmetries.
If the class number $h_K$ of $K$ is equal to one, it is immediate that the two definitions agree, however for $h_K > 1$ their symmetries contain endomorphisms (see Prop. 2.17 \cite{cmr1}) whereas our symmetries are always given by automorphisms. We want to mention that it is no problem to generalize (this is already contained in \cite{conmar}) their definition to the context of a BC-type system for an arbitrary number field and, without changing the definition of our arithmetic subalgebra, we could have proved Theorem \ref{theorem} by using the new definition of symmetries (now containing endomorphisms). \\
This might look odd at first sight but in the context of endomotives (see \cite{conmar} for a reference) the relationship between the two different definitions will become transparent. More precisely in a forthcoming article we will show that every BC-type system (for an arbitrary number field) is an endomotive and the precise relationship between the two definitions of symmetries (and their actions on (extremal) $KMS_\beta$-states) will become clear.
\end{remark}

%%%%%%%%%%%%%%%%%%%%%%%%%%%%%%%%%%%%%%%%%%%%%%%%%%%%
%%%%%%%%%%%%%%%%%%%%%%%%%%%%%%%%%%%%%%%%%%%%%%%%%%%%
%%%%%%%%%%%%%%%%%%%%%%%%%%%%%%%%%%%%%%%%%%%%%%%%%%%%

\appendix

\section{Algebraic Groups}

Our references are \cite{water} and \cite{cm}. Let $k$ denote a field of characteristic zero, $K$ a finite field extension of $k$ and $\overline k$ an algebraic closure of $k$. Further we denote by $R$ a $k$-algebra.

\subsection{Functorial definition and basic constructions}
\label{defA}
An \textbf{(affine) algebraic group} $G$ (over $k$) is a representable functor from (commutative, unital) $k$-algebras to groups. We denote by $k[G]$ its representing algebra, i.e. for any $R$ we have $G(R)=Hom_{k-alg}(k[G],R)$. \\ 
A \textbf{homomorphism} $F: G\to H$ between two algebraic groups $G$ and $H$ (over $k$) is given by a natural transformation of functors. \\ \\
Let $G$ and $H$ be two algebraic groups over $k$, then their \textbf{direct product} $G\times H$ is the algebraic group (over $k$) given by $R \mapsto G(R) \times H(R)$. \\ 
Let $G$ be an algebraic group over $k$ and $K$ and extension of $k$. Then by \textbf{extension of scalars} we obtain an algebraic group $G_K$ over $K$ represented by $K \otimes_k k[G]$.\\
Now let $G$ be an algebraic group over $K$. The \textbf{Weil restriction} $Res_{K/k}(G)$ is an algebraic group over $k$ defined by $Res_{K/k}(G)(R) = G(K \otimes_k R)$.
\begin{remark}
All three constructions are functorial.
\end{remark}

\subsection{Examples}
\label{exA}
1) The \textbf{multiplicative group} $\mathbb G_{m,k}$ (over $k$) is represented by $k[x,x^{-1}]=k[x,y]/ (xy-1)$, i.e. $\mathbb G_{m,k}(R) = R^\times$. \\
2) Define $\mathbb S = Res_{\complex / \real}(\mathbb G _{m,\complex})$. We have $\mathbb S (\real) = \complex ^\times$ and $\mathbb S (\complex) \cong \complex ^\times \times \complex ^\times$. In particular we have $\mathbb S _\complex \cong \mathbb G_{m,\complex}\times \mathbb G_{m,\complex}$. \\
3) More general an algebraic group $T$ over $k$ is called a \textbf{torus} if $T_{\overline k}$ is isomorphic to a product of copies of $\mathbb G_{m,\overline k}$.\\
4) The \textbf{general symplectic group} $GSp(V,\psi)$ attached to a symplectic $\q$-vector space $(V,\psi)$ is an algebraic group over $\q$ defined on a $\q$-algebra $R$ by \begin{align*} &GSp(R) = \\ &\{f \in End_R(V \otimes_\q R) \ | \ \exists \ \nu(f) \in R^\times \text{ : } \psi_{R}(f(x),f(y)) = \nu(f) \psi_{R}(x,y) \ \forall x,y \in V\otimes_\q R\}. \end{align*}

\subsection{Characters}
\label{charA}

Let $G$ be an algebraic group over $k$ and set $\Lambda = \ganz[Gal(\overline k / k)]$. The \textbf{character group} $X^*(G)$ of $G$ is defined by $Hom(G_{\overline k},\mathbb G_{m,\overline k})$. There is a natural action of $Gal(\overline k / k)$ on $X^*(G)$, i.e. $X^*(G)$ is a $\Lambda$-module. Analogously the \textbf{cocharacter group} $X_*(G)$ of $G$ is the $\Lambda$-module $Hom(\mathbb G_{m,\overline k},G_{\overline k})$. We denote the action of $\sigma \in \Lambda$ on a (co)character $f$ by $\sigma f$ or $f^\sigma$. There is the following important 
\begin{theorem}[7.3 \cite{water}]
\label{thmA}
The functor $G \mapsto X^*(G)$ is a contravariant equivalance from the category of algebraic groups of multiplicative type over $k$ and the category of finitely generated abelian groups with a continuous action of $Gal(\overline k / k)$.
\end{theorem} 

\begin{remark}
See 7.2 \cite{water} for the definition of groups of multiplicative type. We only have to know that algebraic tori are of multiplicative type.
\end{remark}
There is a natural bi-additive and $Gal(\overline k / k)$-invariant pairing $<\cdot,\cdot> : X_*(G)\times X^*(G) \to \ganz$ given by $<\chi,\mu> = \mu \circ \chi \in Hom(\mathbb G_{m,\overline k},\mathbb G_{m,\overline k}) \cong \ganz$. If $G$ is of multiplicative type the pairing is perfect, i.e. there is an isomorphism of $\Lambda$-modules $X_*(G) \cong Hom(X^*(G),\ganz)$.

\subsection{Norm maps}
\label{nmA}
Let $L$ be a finite field extension of $K$, i.e. we have a tower $k\subset K \subset L$, and let $T$ be a torus over $k$. Then there are two types of morphisms of algebraic groups which we call \textbf{norm maps}. The first one \eq{Nm_{L/K} : Res_{L/K}(T_L) \to T_K}
is induced by the usual norm map of algebras $R \otimes_K L \to R$, for $R$ a $K$-algebra. In applying the
Weil restriction functor $Res_{K / k}$ we obtain the second one, namely
$$N_{L/K}:Res_{L/k}(T_L)\to Res_{K/k}(T_K).$$

\subsection{The case of number fields}
\label{nfA}
Let $K$ be a number field. We are interested in the algebraic group $T^K=Res_{K/\q}(\mathbb G_{m,K})$ (over $\q$). We have $T^K(R) = (K\otimes_\q R)^\times$. \\ It is easy to see that the isomorphism $K \otimes _\q \overline \q \cong \prod_{\rho \in Hom(K,\overline \q)}\overline \q$ induces an isomorphism of algebraic groups $T^K_{\overline \q} \cong \prod_{\rho \in Hom(K,\overline \q)} \mathbb G_{m,\overline \q}$. It follows that $X^*(T^K) \cong \ganz^{Hom(K,\overline \q)}$ with $Gal(\qbar/ \q)$ acting as follows. \\
For $f = \sum_{\rho \in Hom(K,\qbar)} a_\rho [\rho] \in \ganz^{Hom(K,\qbar)}$ and $\sigma \in Gal(\qbar / \q)$ we have $$\sigma f = \sum_{\rho \in Hom(K,\qbar)} a_\rho [\sigma \circ \rho]=\sum_{\rho \in Hom(K,\qbar)} a_{\sigma^{-1}\circ\rho} [\rho].$$ 
For any inclusion $K \subset L$ of number fields the norm map $N_{L/K} : T^L \to T^K$ is defined by saying that a character $f= \sum_{\rho \in Hom(K,\qbar)} a_\rho [\rho] \in X^*(T^K)$ is mapped to the character $f_L= \sum_{\rho' \in Hom(L,\qbar)} (a_{\rho' | K}) [\rho'] \in X^*(T^L)$.

%%%%%%%%%%%%%%%%%%%%%%%%%%%%%%%%%%%%%%%%%%%%%%%%%%%%
%%%%%%%%%%%%%%%%%%%%%%%%%%%%%%%%%%%%%%%%%%%%%%%%%%%%
%%%%%%%%%%%%%%%%%%%%%%%%%%%%%%%%%%%%%%%%%%%%%%%%%%%%

\section{CM fields}

We follow \cite{cm} and \cite{isv}. By $\iota$ we denote the complex conjugation of $\komplex$.

\subsection{CM fields and CM types}
\label{defB}
Let $E$ denote a number field. If E is a totally imaginary quadratic extension of a totally real field, we call $E$ a \textbf{CM field}. In particular the degree of a CM field is always even. A \textbf{CM type} $(E,\Phi)$ is a CM field $E$ together with a subset $\Phi \subset Hom(E,\complex)$ such that $\Phi \cup \iota \Phi = Hom(E,\complex)$ and $\Phi \cap \iota \Phi =\emptyset$.

\subsection{About $h_\phi$ and $\mu_\phi$}
\label{canmapB}
Let $(E,\Phi)$ be a CM type. Then there are natural isomorphisms $T^E_\real \cong \prod_{\phi \in \Phi} \mathbb S$ resp. $T^E_\complex \cong \prod_{\phi \in \Phi}\mathbb G_{m,\complex} \times \prod_{\overline \phi \in \iota\Phi}\mathbb G_{m,\complex} $, where the first one is induced by $E \otimes_\q \real \cong \prod_{\phi \in \Phi} \komplex$ and the second one by $E \otimes \complex \cong \prod_{\phi \in \Phi} \complex \times \prod_{\overline \phi \in \iota \Phi} \complex$. \\
Thus we obtain natural morphisms \eq{ \label{hphiB} h_\Phi : \mathbb S \to T^E_\real \text{ ; $z \mapsto (z)_{\phi \in \Phi}$} } and \eq{ \label{muphiB} \mu_\Phi : \mathbb G _{m,\complex} \to T^E_\complex \text{ ; $z \mapsto (z)_{\phi \in \Phi} \times (1)_{\overline \phi \in \iota \Phi} $.} }
If we take $\mu_\Phi$ for granted we could have defined $h_{\Phi}$ by the composition \diagram{Res_{\complex / \real}(\mathbb G_{m,\complex}) \ar[rr]^{Res_{\complex / \real}(\mu_{\Phi})} & & Res_{\complex / \real} (T^E_\complex) \ar[rr]^{\ \ \ \ Nm_{\complex / \real}} & & T^E_\real.}
In particular we see that $h_\Phi$ and $\mu_\Phi$ are related by \eq{h_{\Phi,\complex}(z,1)=\mu_\Phi(z).}

\begin{remark}
In the last two sections one might have replaced $\complex$ by $\qbar$.
\end{remark}

\subsection{The reflex field and reflex norm}
\label{reflexB}
Let $(E,\Phi)$ be a CM type. The \textbf{reflex field} $E^*$ of $(E,\Phi)$ is the subfield of $\overline {\q}$ defined by any one of the following conditions: \\
a) $\sigma \in Gal(\overline \q / \q)$ fixes $E^*$ if and only if $\sigma \Phi = \Phi$ ; \\
b) $E^*$ is the field generated over $\q$ by the elements $\sum_{\phi\in \Phi}\phi(e)$, $e \in E$ ; \\
c) $E^*$ is the smallest subfield of $\overline \q$ such that there exists a $E\otimes_\q E^*$-module $V$ such that \eq{Tr_{E^*}(e | V) = \sum_{\phi \in \Phi}\phi(e) \text{, for all $e \in E.$}}
The \textbf{reflex norm} of $(E,\Phi)$ is the morphism of algebraic groups $N_{\Phi}:T^{E^*}\to T^E$ given, for $R$ a $\q$-algebra, by \eq{ \label{reflexexplicit} a \in T^{E^*}(R) \mapsto det_{E\otimes_\q R}(a | V\otimes_\q R) \in T^E(R). }

%%%%%%%%%%%%%%%%%%%%%%%%%%%%%%%%%%%%%%%%%%%%%%%%%%%%
%%%%%%%%%%%%%%%%%%%%%%%%%%%%%%%%%%%%%%%%%%%%%%%%%%%%
%%%%%%%%%%%%%%%%%%%%%%%%%%%%%%%%%%%%%%%%%%%%%%%%%%%%

\section{The Serre group}
\label{serre group}

Our references are \cite{pedestrian}, \cite{cm} and \cite{wei}. Let $K$ be a number field. We fix an embedding $\tau : K \to \qbar \to \complex$ and denote by $\iota$ complex conjugation on $\complex$.

\subsection{Definition of the Serre group}
\label{defC}
The following are equivalent: \\ \\
(1) The \textbf{Serre group} attached to $K$ is a pair $(S^K,\mu^K)$ consisting of a $\Q$-algebraic torus $S^K$ and a cocharacter $\mu^K \in X_*(S^K)$ defined by the following universal property.
For every pair $(T,\mu)$ consisting of a $\Q$-algebraic torus $T$ and a cocharacter $\mu \in X_*(T)$ defined over $K$ satisfying the Serre condition \eq{\label{serrecond}(\iota + 1)(\sigma-1)\mu =0=(\sigma-1)(\iota+1)\mu \ \ \forall \sigma \in Gal(\overline \Q / \Q)}
there exists a unique morphism $\rho_\mu : S^K \to T$ such that the diagram
\diagramst{S^K_\qbar \ar[rr]^{\rho_{\mu,\qbar}} & & T_\qbar \\ & \mathbbm G _{m,\qbar} \ar[ul]^{\mu^K} \ar[ur]^{\mu} & }
commutes. \\ \\
(2) The \textbf{Serre group} $S^K$ is defined to be the quotient of $T^K$ such that $X^*(S^K)$ is the subgroup of $X^*(T^K)$ given by all elements $f \in X^*(T^K)$ which satisfy the Serre condition \eqst{(\sigma -1)(\iota+1)f = 0 = (\iota+1)(\sigma-1)f \ \ \forall \sigma \in Gal(\qbar/\q).}
The cocharacter $\mu^K$ is induced by the cocharacter $\mu_\tau \in X_*(T^K)$ defined by \eqst{<\mu_\tau, \Sigma n_\sigma [\sigma]> = n_{\tau}, \ \ \forall \ \Sigma n_\sigma [\sigma] \in \ganz^{Hom(K,\qbar)} \cong X^*(T^K).}
(3) If $K$ does not contain a CM subfield, we set $E= \q$, otherwise $E$ denotes the maximal CM subfield of $K$ and $F$ the maximal totally real subfield of $E$. Then there is an exact sequence of $\q$-algebraic groups \diagram{1 \ar[r] & ker(N_{F/\q} : T^F \to T^\q)\ar[r]^{ \ \ \ \ \ \ \ \ \ \ \ \ i}&T^K \ar[r]^{\pi^K} & S^K \ar[r] & 1,} where $i$ is the obvious inclusion. The cocharacter $\mu_\tau$ of $T^K$, defined as in (2), induces $\mu^K$, i.e. $ \mu^K=\pi^K \circ\mu_\tau$.
\begin{remark}
For $K = \q$ or $K$ an imaginary quadratic field there is the obvious equality $S^K = T^K$.
\end{remark}

\subsection{About $\mu^K$ and $h^K$}
\label{canmapC}
The cocharacter $\mu^K= \pi^K \circ \mu_\tau : \mathbb G _{m,\complex} \to S^K_\complex$ from the last section induces a natural morphism \eq{ \label{hKC} h^K : \mathbbm S \to S^K_\reell } defined by \diagram{Res_{\komplex / \real}(\mathbb G _{m,\komplex}) \ar[rr]^{\ \ Res_{\complex / \real}(\mu^K) \ \ } & & Res_{\komplex / \real}(S^K_\komplex) \ar[rr]^{ \  \ \ \  Nm_{\komplex / \real}} & & S^K_\real. }  
We see that $\mu^K$ and $h^K$ are related by \eq{h^K_\complex (z,1) = \mu^K(z)}
or in other words, for $z \in \complex^\times$, we have $h^K(z) = \mu^K(z)\mu^K(z)^\iota.$

\subsection{About $\rho_\Phi$ and the reflex norm $N_\Phi$}
\label{phirho}
Let $(E,\Phi)$ be a CM type. The natural morphism $\mu_\Phi \in X_*(T^E)$ (cf. \ref{muphiB}) is defined over the reflex field $E^*$ and an easy calculation shows that it satisfies the Serre condition (\ref{serrecond}). By the universal property of the Serre group we obtain a $\q$-rational morphism \eq{\rho_\Phi : S^{E^*} \longrightarrow T^E} such that \eq{\mu_\Phi = \rho_{\Phi,\complex} \circ \mu^{E^*}.} Also, we see immediately that \eq{\label{comphphi}h_\Phi = \rho_{\Phi,\real} \circ h^{E^*}.}
Moreover we can relate $\rho_\Phi$ and $\mu_\Phi$ by the following commutative diagram \diagram{T^{E^*}\ar[rr]^{Res(\mu_{\Phi,E^*}) \ \ \ \ \ } \ar[dr]^{\pi^{E^*}} & & Res_{E^*/\q}(T^E_{E^*}) \ar[rr]^{\ \ \ \ \ \ \ N_{E^*/\q}} & & T^E \\ & S^{E^*} \ar[rrru]_{\rho_\Phi}} 
which can be seen on the level of characters.
The relation with the reflex norm $N_\Phi : T^{E^*} \to T^E$ (see \ref{reflexB}) is given by the following important

\begin{proposition}[\cite{cm}]
\label{propreflexnorm}
We have the equality \eq{\label{reflexnormprop}N_\Phi = N_{E^* / \q} \circ Res(\mu_{\Phi,E^*}).}
\end{proposition}

\subsection{More properties of the Serre group}
\label{propC}
The following properties are all taken from \cite{cm}.
\begin{proposition}
\label{propC1}
Let $E \subset K$ denote two number fields.
\begin{enumerate}
\item The norm map $N_{K/E}:T^K\to T^E$ induces a commutative diagram
\diagram{T^K \ar[r]^{N_{K/E}} \ar[d]^{\pi^K} & T^E \ar[d]^{\pi^E} \\ S^K \ar[r]^{} & S^E.}
We call the induced morphism $N_{K/E} : S^K \to S^E$.
\item There is a commutative diagram \diagramlab{comphnorm}{\mathbb S \ar[r]^{h^K} \ar[dr]^{h^E} & S^K \ar[d]^{N_{K/E}} \\ & S^E}

\item Let $E$ denote the maximal CM field contained in $K$, if there is no such subfield we set $E=\q$. Then $N_{K/E} : S^K \to S^E$ is an isomorphism.

\item Let $(E,\Phi)$ be a CM type and $K_1 \subset K_2$ two number fields, such that $E^* \subset K_1$, and let $\rho_{\Phi,i} : S^{K_i} \to T^E$ be the corresponding maps from the universal property of the Serre group. Then we have $$\rho_{\Phi,1} \circ N_{K_2 / K_1} = \rho_{\Phi,2}.$$
\end{enumerate}
\end{proposition}

%%%%%%%%%%%%%%%%%%%%%%%%%%%%%%%%%%%%%%%%%%%%%%%%%%%%
%%%%%%%%%%%%%%%%%%%%%%%%%%%%%%%%%%%%%%%%%%%%%%%%%%%%
%%%%%%%%%%%%%%%%%%%%%%%%%%%%%%%%%%%%%%%%%%%%%%%%%%%%

\section{Shimura Varieties}
Our references are Deligne's foundational \cite{Del2}, Milne's \cite{isv} and Hida's \cite{hida}. \\ \\
Let $G$ be an algebraic group over $\q$. Then the adjoint group $G^{ad}$ of $G$ is defined to be the quotient of $G$ by its center $C$. The derived group $G^{der}$ of $G$ is defined to be the intersection of the normal algebraic subgroups of $G$ such that $G / N$ is commutative. By $G(\real)^+$ we denote the identity component of $G(\real)$ relative to its real topology and set $G(\q)^+ = G(\q) \cap G(\real)^+$. If $G$ is reductive, we denote by $G(\real)_+$ the group of elements of $G(\real)$ whose image in $G^{ad}(\real)$ lies in its identity component and set $G(\q)_+ = G(\q) \cap G(\real)_+$.

\subsection{Shimura datum}
\label{datumD}
A \textbf{Shimura datum} is a pair $(G,X)$ consisting of a reductive group $G$ (over $\q$) and a $G(\real)$-conjugacy class $X$ of homomorphisms $h:\mathbb S \to G_\real$, such that the following (three) axioms are satisfied \\ \\
$(SV1)$: For each $h \in X$, the representation $Lie(G_\real)$ defined by $h$ is of type $\{(-1,1),(0,0),(1,-1)\}$. \\
$(SV2)$: For each $h \in X$, $ad(h(i))$ is a Cartan involution on $G_\real^{ad}$. \\
$(SV3)$: $G^{ad}$ has no $\q$-factors on which the projection of $h$ is trivial. \\ \\
Because $G(\real)$ is acting transitively on $X$ it is enough to give a morphism $h_0 : \mathbb S \to G_\real$ to specify a Shimura datum.
Therefore a \textbf{Shimura datum} is sometimes written as triple $(G,X,h_0)$ or simply by $(G,h_0)$. \\ 
Further in our case of interest the following axioms are satisfied and (simplify the situation enormously). \\ \\
$(SV4)$: The weight homomorphism $\omega_X : \mathbb G_{m,\real} \to G_\real$ is defined over $\q$. \\
$(SV5)$: The group $C(\q)$ is discrete in $C(\adfin)$. \\
$(SV6)$: The identity component of the center $C^o$ splits over a CM-field. \\
$(SC)$: The derived group $G^{der}$ is simply connected. \\
$(CT)$: The center $C$ is a cohomologically trivial torus. \\
\begin{remark}
1) Axioms $(SV1-6)$ are taken from \cite{isv}, the other two axioms are taken from \cite{hida}. \\
2) The axioms of a Shimura variety ($SV1$-$3$) imply, for example, that $X$ is a finite union of hermitian symmetric domains. When viewed as an analytic space we sometimes write $x$ instead of $h$ for points in $X$ and $h_x$ for the associated morphism $h_x : \mathbb S \to G_\real$. \\
3) In 3.1 \cite{HaPa} a more general definition of a Shimura datum is given. For our purpose Deligne's original definition, as given above, and so called $0$-dimensional Shimura varieties are sufficient. \\
\end{remark}
A \textbf{morphism of Shimura data} $(G,X) \to (G',X')$ is a morphism $G\to G'$ of algebraic groups which induces a map $X \to X'$.

\subsection{Shimura varieties}
\label{shvD}
Let $(G,X)$ be a Shimura datum and let $K$ be a compact open subgroup of $G(\adfin)$. Set $Sh_K=Sh_K(G,X) = G(\q) \backslash X \times G(\adfin) / K,$
where $G(\q)$ is acting on $X$ and $G(\adfin)$ on the left, and $K$ is acting on $G(\adfin)$ on the right. On can show (see 5.13 \cite{isv}) that there is a homeomorphism $Sh_K \cong \bigsqcup \Gamma_g \backslash X^+$.
Here $X^+$ is a connected component of $X$ and $\Gamma_g$ is the subgroup $gKg^{-1} \cap G(\q)_+$ where $g$ runs runs over a set of representatives of $G(\q)_+ \backslash G(\adfin) / K$. When $K$ is chosen sufficiently small, then $\Gamma_g \backslash X^+$ is an arithmetic locally symmetric variety. For an inclusion $K' \subset K$ we obtain a natural map $Sh_{K'} \to Sh_K$ and in this way an inverse system $(Sh_K)_K$.
There is a natural right action of $G(\adfin)$ on this system (cf. p 55 \cite{isv}).  \\
The \textbf{Shimura variety} $Sh(G,X)$ associated with the Shimura datum $(G,X)$ is defined to be the inverse limit of varieties $\varprojlim_K Sh_K(G,X)$ together with the natural action of $G(\adfin)$. Here $K$ runs through sufficiently small compact open subgroups of $G(\adfin)$. $Sh(G,X)$ can be regarded as a scheme over $\complex$.\\
Let $(G,X)$ be a Shimura datum such that $(SV5)$ holds, then one has \eq{\label{sv5}Sh(G,X) = \varprojlim _K Sh_K(G,X) = G(\q)\backslash X \times G(\adfin).}
In this case we write $[x,l]$ for an element in $Sh(G,X)$ and the (right) action of an element $g \in G(\adfin)$ is given by \eq{g[x,l] = [x,lg].} In the general case, wenn $(SV5)$ is not holding we use the same notation, understanding that $[x,l]$ stands for a family $(x_K,l_K)_K$ indexed by compact open subgroups $K$ of $G(\adfin)$. \\
A \textbf{morphism of Shimura varieties} $Sh(G,X)\to Sh(G',X')$ is an inverse system of regular maps of algebraic varieties compatible with the action of $G(\adfin)$. We have the following functorial property:\\
A morphism $\varphi : (G,X)\to (G',X')$ of Shimura data defines an equivariant morphism $Sh(\varphi) : Sh(G,X)\to Sh(G',X')$ of Shimura varieties, which is a closed immersion if $G\to G'$ is injective (Thm 5.16 \cite{isv}).

\subsection{Example}
\label{exD}
We want to give some details about the Shimura varieties attached to the data $\mathcal S_\siegel$ constructed in \ref{construction}. For the identification of the $GSp(\real)$-conjugacy class of $h_{cm}$ with the higher Siegel upper lower half space \alst{\mathbb H_g^±=\{ M = A + i B \in M_g(\complex) \ | \ A = A^t, B \text{ positive or negative definitive } \} }
we refer further to exercise 6.2 \cite{isv}. \\
In addition the data $\mathcal S_\siegel$ fulfill all the axioms stated in \ref{datumD}. The validity of $(SV1$-$6)$ is shown on p. 67 \cite{isv} and the validity of $(SC)$ and $(CT)$ in \cite{hida}. The latter two axioms are important for making the arguments in \cite{MiS} in this case. 
From $(SV5)$ follows in particular that we don't have to bother about the limits in the definition of $Sh(GSp,\mathbb H_g^±)$ because we have \alst{Sh(GSp,\mathbb H_g^±) = GSp(\q) \backslash (\mathbb H _g^± \times GSp(\adfin)).}

\subsection{Connected Shimura varieties}
\label{cshvD}
A \textbf{connected Shimura datum} is a pair $(G,X^+)$ consisting of a semisimple algebraic group $G$ over $\q$ and a $G^{ad}(\real)^+$-conjugacy class of homomorphisms $h:\mathbb S \to G^{ad}_\real$ satisfying axioms 
$(SV1$-$3)$. \\
The \textbf{connected Shimura variety} $Sh^o(G,X^+)$ associated with a connected Shimura datum $(G,X^+)$ is defined by the inverse limit \begin{align} Sh^o(G,X^+) = \varprojlim_\Gamma \Gamma \backslash X^+ \end{align} where
$\Gamma$ runs over the torsion-free arithmetic subgroups of $G^{ad}(\q)^+$ whose inverse image in $G(\q)^+$ is a congruence subgroup. \\ \\
If we start with a Shimura datum $(G,X)$ and choose a connected component $X^+$ of $X$, we can view $X^+$ as a $G^{ad}(\real)^+$-conjugacy class of morphisms $h : \mathbb S \to G^{ad}_\real$ by projecting elements in $X^+$ to $G^{ad}_\real$. One can show that $(G^{der},X^+)$ is a connected Shimura datum.
Further if we choose the connected component $Sh(G,X)^o$ of $Sh(G,X)$ containing $X^+ \times {1}$, one has the following compatibility relation \begin{align}Sh(G,X)^o = Sh^o(G^{der},X^+).\end{align}

\subsection{$0$-dimensional Shimura varieties}
\label{0shvD}
In section \ref{shdatum1} we defined a "Shimura datum" $\mathcal S_K = (T^K,X_K)$ which is not a Shimura datum in the above sense because $X_K$ has more than one conjugacy class (recall that $T^K$ is commutative). Rather $\mathcal S_K$ is a Shimura datum in the generalized sense of Pink \cite{pink} which we don't want to recall here. Instead we define the notion of a $0$-dimensional Shimura varieties following \cite{isv} which covers all exceptional Shimura data we consider. We define a \textbf{$0$-dimensional Shimura datum} to be a triple $(T,Y,h)$, where $T$ is a torus over $\q$, $Y$ a finite set on which $T(\real) / T(\real)^+$ acts transitively and $h : \mathbb S \to T_\real$ a morphism of algebraic groups. We view $Y$ as a finite cover of $\{h\}$. We remark that the axioms ($SV1-3$) are automatically satisfied in this setup.\\
The associated \textbf{$0$-dimensional Shimura variety} $Sh(T,Y,h)$ is defined to be the inverse system of finite sets $T(\q) \backslash Y \times T(\adfin) / K$ with $K$ running over the compact open subgroups of $T(\adfin)$. \\
A \textbf{morphism} $(T,Y,h) \to (H,h_0)$ from a $0$-dimensional Shimura datum to a Shimura datum, with $H$ an algebraic torus, is given by a morphism of algebraic groups $\varphi : T \to H$ such that $h = \varphi_\real \circ h_0$. \\
If $\varphi$ is such a morphism it defines a morphism $Sh(\varphi) : Sh(T,Y,h) \to Sh(H,h_0)$ of Shimura varieties.

\begin{remark}
\label{easywritingsk}
We have that $\mathcal S _K$ fulfills axiom $(SV5)$ if and only if $K=\q$ or $K$ an imaginary quadratic field (see 3.2 \cite{HaPa}). 
\end{remark}

\subsection{Canonical model of Shimura varieties}
\label{canmodD}
Let $(G,X)$ be a Shimura datum. A point $x \in X$ is called a \textbf{special point} if there exists a torus $T \subset G$ sucht that $h_x$ factors through $T_\real$. The pair $(T,x)$ or $(T,h_x)$ is called special pair. If $(G,X)$ satisfies the axioms $(SV4)$ and $(SV6)$, then a special point is called \textbf{CM point} and a special pair is called CM pair. \\ \\
Now given a special pair $(x,T)$ we can consider the cocharacter $\mu_x$ of $G_\complex$ defined by $\mu_x(z) = h_{x,\complex}(z,1)$. Denote by $E(x)$ the field of definition of $\mu_x$, i.e.  the smallest subfield $k$ of $\complex$ such that $\mu_x : \mathbb G_{m,k} \to G_k$ is defined. \\
Let $R_x$ denote the composition \diagram{T^{E(x)} \ar[rr]^{Res_{E(x) / \q}(\mu_x) \ \  \ \ \ \ \ \ } & & Res_{E(x)/\q}(T_E(x)) \ar[rr]^{\ \ \ \ \ \ \ \ \ \ Nm_{E(x)/ \q}} & & T } and define the \textbf{reciprocity morphism} \eq{r_x = R_x(\adfin) : \mathbb A ^\times_{E(x),f} \to T(\adfin).}
Moreover every datum $(G,X)$ defines an algebraic number field $E(G,X)$, the reflex field of $(G,X)$. For the definition we refer the reader to 12.2 \cite{isv}. 
\begin{remark}
\label{remcanmodD}
1) For the Shimura datum $\mathcal S _\siegel = (GSp,\mathbb H_g^±)$ (see \ref{construction}) we have $E(\mathcal S_K) = \q$ (cf. p. 112 \cite{isv}). \\
2) For explanations to relations with the reflex field of a CM field (cf. \ref{reflexB}), see example 12.4 b) of pp. 105 \cite{isv}.
\end{remark}
A \textbf{model} $M^o(G,X)$ of $Sh(G,X)$ over the reflex field $E(G,X)$ is called \textbf{canonical} if \\
1)$M^o(G,X)$ is equipped with a right action of $G(\adfin)$ that induces an equivariant isomorphism $M^o(G,X)_\complex \cong Sh(G,X)$, and \\
2) for every special pair $(T,x) \subset (G,X)$ and $g \in G(\adfin)$ the point $[x,g] \in M^o(G,X)$ is rational over
$E(x)^{ab}$ and the action of $\sigma \in Gal(E(x)^{ab} / E(x))$ is given by \begin{align} \label{canmodspD} \sigma [x,g] = [x, r_x(\nu)g] \end{align} where $\nu \in \mathbb A^\times _{E(x),f}$ is such that $[\nu] = \sigma^{-1}$ under Artin's reciprocity map. \\
In particular, for every compact open subgroup $K \subset G(\adfin)$, it follows that $M_K^o(G,X)=M^o(G,X) / K$ is a model of $Sh_K(G,X)$ over $E(G,X)$.
\begin{remark}
Canonical models are known to exist for all Shimura varieties (see \cite{isv}).
\end{remark}
\subsection{Canonical model of connected Shimura varieties}
\label{ccanmodD}
We refer to 2.7.10 \cite{Del2} for the precise definition of the canonical model $M^o(G,X^+)$ of a connected Shimura variety $Sh(G,X^+)$. Here we just want to mention the compatibility \begin{align}M^o(G^{der},X^+)=M^o(G,X)^o\end{align}
where the latter denotes a correctly chosen connected component of the canonical model $M^o(G,X)$.

\bibliographystyle{alpha} 
\bibliography{bcandcm}

\end{document}